\documentclass[a4paper]{amsart}
\usepackage[english]{babel}
\usepackage{amsmath, amsthm, amsfonts, amssymb, mathtools}
\usepackage{accents}
\usepackage[final]{hyperref}

\newtheorem{theorem}{Theorem}

\theoremstyle{remark}
\newtheorem{remark}{Remark}
\newtheorem{example}{Example}

\numberwithin{equation}{section}

\allowdisplaybreaks
\newcommand{\e}{\mathrm{e}}
\newcommand{\1}{\!\!\,{\textbf 1}}
\newcommand{\ov}{\overline}

\newcommand{\wt}{\widetilde}
\newcommand{\wh}{\widehat}
\newcommand{\pt}{\partial}
\newcommand{\Id}{\mathrm{Id}}

\newcommand{\Var}{\mathrm{Var}}
\newcommand{\p}{\partial}

\DeclareMathOperator\E{E}

\DeclareMathOperator\Prob{P}

\DeclareMathOperator{\grad}{grad}
\DeclareMathOperator{\Div}{div}
\DeclareMathOperator{\Hess}{Hess}
\DeclareMathOperator{\tr}{tr}

\newcommand{\ve}{\varepsilon}

\newcommand{\mbR}{{\mathbb R}}

\newcommand{\mbN}{{\mathbb N}}
\newcommand{\mbC}{{\mathbb C}}

\newcommand{\dt}[1]{\accentset{\approx}{#1}}

\makeatletter
\@namedef{subjclassname@2020}{
  \textup{2020} Mathematics Subject Classification}
\makeatother

\title[An introduction to operator splitting methods]{A quick probability-oriented introduction to operator splitting methods}
\author{M.B. Vovchanskyi}
\address{Institute of Mathematics, National Academy of Sciences of Ukraine, Teresh\-chenkivska Str. 3, Kyiv 01601, Ukraine}
\email{vovchansky.m@gmail.com}


\usepackage[backend=biber, 
	bibencoding=utf8,
	style=ieee, 
	citestyle=numeric-comp,
	sorting=nty,
	autolang=other,
	maxcitenames=4,
	maxbibnames=10,
	giveninits=true,
	date=year,	doi=false,isbn=false,url=false,
	eprint=true,
	hyperref=true
	]{biblatex}
\AtEveryBibitem{\clearfield{language}}
\AtEveryBibitem{\clearfield{note}}
 
\usepackage{hyperref}

\begin{document} 

\begin{abstract}
The survey is devoted to operator splitting methods in the abstract formulation and their applications in probability. While the survey is focused on multiplicative methods, the BCH formula is used to discuss exponential splitting methods and a short informal introduction to additive splitting is presented. We introduce frameworks and available deterministic and probabilistic results and concentrate on constructing  a wide picture of the field of operator splitting methods, providing a rigorous description in the setting of abstract Cauchy problems and an informal discussion for further and parallel advances. Some limitations and common difficulties are listed, as well as examples of works that provide solutions or hints. No new results are provided. The bibliography contains illustrative deterministic examples and a selection of probability-related works. 
\end{abstract}

\maketitle 

\tableofcontents

\section{The scope  of the paper and first examples}

{\em This paper is an extended and reworked version of a short course given by the author at ''Uzbekistan-Ukrainian readings in stochastic processes'', Tashkent-Kyiv, 2022, and was prepared for a special issue of ''Theory of stochastic processes'', devoted to publishing lecture notes from the aforementioned workshop.}

Operator splitting methods is a family of well-known methods of decomposing a dynamical system by providing a representation of the governing mechanism as a sum of simpler components (''forces'') and using this representation to provide an approximation of the real trajectory of the dynamical system.  Though this idea of ''splitting'' is often associated with the celebrated abstract Trotter-Kato formula, in which case the approximation is constructed by using  a composition of subsystems each of which is driven by exactly one component of the aforementioned representation, it goes far beyond that and includes weighted (linear) combinations of such subsystems, compositions of mixed backward-forward Euler schemes and can be encountered in the optimization theory in a disguise. 

More specifically, the idea is as follows. The time interval is divided into sufficiently small steps, and  on every step a  solution operator  of the original system is replaced with an  approximation that is constructed by a splitting procedure, after which the solutions are combined recursively starting from time $0.$ 

Due to the general nature of the idea, operator splitting can essentially be used anywhere where ODEs or (S)PDEs with a natural decomposition arise, be this decomposition dictated by the presence of co-existing physical, chemical or biological mechanisms, often acting on different space/time scales,  or by properties of available numerical methods. 

\begin{example}
Consider a general advection-diffusion-reaction equation\footnote{Alternative names are general scalar transport equations, convection-diffusion-reaction or convection-diffusion equations etc.  The  term ''advection'' is often used synonymously to ''convection''.}:
\begin{equation}
\label{eq:0.1}
\frac{\partial u}{\partial t} = \Div(A  \grad u) - \Div (Bu) + R,
\end{equation}
where $A$ represents diffusion, the field $B$ is (possibly superficial) velocity and $R$ is sources and sinks. $R$ may be of chemical origin, while the velocity field  reflects physical properties of a reservoir etc.  One may prefer drastically  different numerical methods  for integrating subsystems obtaining by splitting the RHS into a sum of  differential operators of at most second order and treating those subsystems separately. For instance, physically justified first order equations have a developed theory that often explicitly relies on conservation laws (e.g.~\cite{HundsVer07Numerical}) while the diffusion part is typically treated via finite elements (e.g.~\cite{ZienTayZhu13finite}). 
\end{example}

\begin{example}\label{ex:a.2} Assume that $a_1, a_3<0, a_2 \in\mbR, |a_1| \gg |a_2|+|a_3|$ and $|a_2|, |a_3| \approx 1.$ Consider the ODE 
\begin{align*}
& \frac{dy}{dt} = A y, \\
& A = \begin{pmatrix}
a_1 & a_2 \\
a_2 & a_3
\end{pmatrix}.
\end{align*} 
If we try to simulate $y$ by using the explicit forward Euler method with step size $h,$ a necessary condition is 
\[
\|\mathrm{I} + hA\|<1
\]
in the  operator matrix norm, which implies 
that $h$ should be of order $\frac{1}{|a_1|}.$ Since $h$ is limited in practice, it may not be possible for the forward Euler method to be stable. 
However, for the decomposition
\[
A=A_1+A_2 = \begin{pmatrix}a_1 & 0 \\ 0 & 0 \end{pmatrix} + \begin{pmatrix}0 & a_2 \\ a_2 & a_3 \end{pmatrix}
\]
both equations
\[
\frac{dy^{(k)}}{dt} =A_k y^{(k)}, \quad k=1,2,
\]
can be solved explicitly.  This is a toy model example where a stiff ($A_1$) and a non-stiff part ($A_2$) of the original equation can be isolated. Here stiffness (e.g.~\cite{HaiWan96solving, Lam91numerical, HundsVer07Numerical} in the deterministic setting and \cite{Mil95numerical} in the stochastic setting) in the context of operator splitting methods is understood in a vague general sense: a system is stiff if a successful application of certain  numerical methods (usually explicit ones and/or those used in the field as standard integrators) requires an impractically small time step and is thus infeasible or  costly/hard to implement. The idea of the current example combines well  with the previous example since competing mechanisms of different origin (chemical vs. mechanical) may act on different  scales, which is a typical source of stiffness as it directly increases  the stiffness ratio. 
\end{example}

\begin{example}(e.g.~\cite{BlaCa16concise})
As a development of Example~\ref{ex:a.2}, consider a stiff ODE in $\mbR^d$ with constant coefficients
\[
\frac{dy}{dt} = A y,
\] 
where eigenvalues $a_k, k =\ov{1,d},$ of $A$ are distinct, real and negative and  $|a_1| \gg \sum_{k=\ov{2,d}} |a_k|,$ $a_d \gg a_k, k=\ov{1,d-1}.$ Then 
\[
y(t) \sim C \e^{a_d t} \to 0, \quad t\to\infty, 
\]
for any $y(0)$ and the term for $a_1$  gets negligible in comparison to $\e^{a_d t}$ as $t$ grows. But the explicit Euler method with step size $h$ displays the same asymptotic behavior only if $h |a_1| < 2$\footnote{Unless the initial conditions is $0$ in the direction of the eigenvector of $a_1.$}. Moreover, with this condition violated, the approximations diverge as $t\to\infty.$  Thus the stability of the explicit Euler method depends on parameters of a subsystem whose contribution to the dynamics of the whole system is minor at best.  This is a typical feature of stiff systems. Separation of the problematic direction from the rest of the system is beneficial in this case.
\end{example}

\begin{example}
\label{ex:a.3}
Let $a_1, a_2 \in \mbR$. Consider 
\begin{align*}
\frac{\pt u(x,t)}{\pt t} &= \frac{1}{2} \left(a_1\frac{\pt^2}{\pt x_1^2} + a_2\frac{\pt^2}{\pt x_2^2} \right) u(x,t), \quad x \in \mbR^2, t\ge 0, \\
u(x,t) &= u_0(x).
\end{align*}
Let $g_t$ be the one-dimensional Gaussian density with variance $t.$ Then 
\begin{align*}
T^{(k)}_t f(z) &= \int_{\mbR} g_{t / \sqrt{a_k}}(z - x_0) v_0(x_0) dx_0, 
\end{align*}
solve 
\[
\frac{\pt v(z,t)}{\pt t} = \frac{a_k}{2} \frac{\pt^2}{\pt z^2} v(z,t), \quad z \in\mbR, t \ge 0, k=1,2.
\]
Combining operators $T^{(1)}_t$ and $T^{(2)}_t$ on the interval $[0;n^{-1}]$ gives for any $x=(x_1, x_2)$
\[
T^{(1)}_{t/n} \Big( T^{(2)}_{t/n} u_0(\cdot, x_2) \Big) = \E u_0 \Big( A w\Big(\frac{t}{n}\Big) + x\Big), 
\]
where $w$ is a two-dimensional Wiener process and
\[
A = \begin{pmatrix}
a_1 & 0 \\
0 & a_2
\end{pmatrix}.
\]
Thus for any $m\in\mbN$\footnote{$f\circ g = f(g)$ hereinafter.}
\[
\Big( T^{(1)}_{t/n} \circ T^{(2)}_{t/n}\Big)^{n} u_0 = \E u_0 \left( Aw(t) + x\right) = u(x,t),
\]
and the operator splitting scheme is exact in this trivial case. This is an example of the alternating direction method whose name is quite self-explanatory.
\end{example}

As a result, operator splitting procedures are often developed for equations and problems in a general formulation (for instance, for the Navier-Stokes equation, the Zakai equation of the filtration theory, advection-diffusion-reaction equations, and composite optimization problems as an example in a non-PDE  setting).     

Such wide availability of the methods results in applications in physics (mechanics of fluids, gases and solids; classical, quantum and celestial mechanics; electromagnetism; symplectic integration etc.), chemistry, biology, geology, ecology, weather forecasting, finances, general machine learning problems (including image processing, large scale optimization etc.)  and other numerous fields. At the same time, operator splitting methods can be encountered in purely theoretical studies. Proper examples will be given later in the text.

The scope of this survey is purely pedagogical:
\begin{itemize}
\item to suggest a brief and quick introduction to operator splitting methods, presenting basics in a form accessible  for a  newcomer (with a background in probability) while still discussing some limitations and technical issues of such methods;
\item to provide a wider picture and to at least mention directions the whole theory goes beyond the Trotter-Kato formula;
\item to give examples of results developed through the use of operator splitting methods in the field of probability and to show that randomness leads to new insights and new techniques; 
\end{itemize} 
This ''probabilistic'' orientation simply means that we prefer probabilistic examples, assume a standard level of knowledge in probability, concentrate on multiplicative/expo\-nen\-ti\-al splitting 
and include a short survey of the general theory of multiplicative splitting for SPDEs. 

At the same, this also leads to some asymmetries in the text: the discussion of the  Baker-Campbell-Hausdorff formula is rather lengthy only because such material is not expected to be covered in textbooks on probability while the reader is assumed to be familiar with SPDEs and variational methods for them so the part devoted to SPDEs does not explain the setting. 

It should be emphasized that there is more than enough information around for a person interested in the topic, including high quality introductory level texts and textbooks. 
 {\bf Not even trying to compete with classical texts,\footnote{There are also highly efficient surveys (e.g.~\cite{CsoSi21numerical,NeidSteZa18operator,BlaCaMu08splitting}).} we merely concentrate on gathering (at least in the form of bibliographic references) important basic facts, hints and probabilistic results as expanded lecture notes (for the aforementioned course) and  giving a total newcomer basic understanding and a map to navigate. }

The level of exposition is thus more or less basic and shallow even though the results cited are often deep and extremely technical to obtain: given the scope and goals of the survey, it would be irreparably pretentious, futile and plainly impossible to proceed otherwise. 
 Thus we refer to original sources for proofs, details and precise formulations.

 All results discussed in the paper are known.




Obviously, the sheer amount of publications devoted to theoretical and applied studies of operator splitting methods  immediately renders impossible the task of providing an exhaustive or even only a representative bibliography unless attention is restricted to a limited partial case or a specific problem -- and the task is highly nontrivial even then. As a result, the choice of illustrative examples 
is to some extend random, and we apologize in advance to all those whose contribution has not been represented 
and  accept all the critique for such gaps, claiming no ill intent. We also give links to surveys and collection of references instead of citing sources directly on some occasions. Still, we hope that the bibliography is of some independent value. No real attempt to give historical notes or to track results to initial sources has been made.

A reader interested in a general and broad picture may be proposed to consult the following selections of works as entry points: 
\begin{itemize}
\item \cite{Pazy12Semigroups, Goldstein85semigroups, EnNa00One-parameter, Da0One-parameter, ItoKa02evolution, Cher74product, EthKur09Markov} for the Trotter-Kato formula, the perturbation theory for semigroups and the formulation in the terms of abstract evolution equations,
\item \cite{GlowOshYin17Splitting, FaHa09splitting}\footnote{Chapters 1--2 for a general exposition and a survey.} for a survey, the general theory, history, connections to optimization problems, PDEs and variational inequalities, and examples,
\item \cite{HundsVer07Numerical} for a treatment of advection-diffusion-reaction equations, examples and a discussion of numerics,
\item  \cite{ReedSai72methods, Si05functional, JohnLa00Feynman, LoHiBetz20FeynmanKac, Hall13quantum, ReedSi75methods_2, Zeid95applied} for the physical exposition and applications to Feynman path integrals and the Feynman-Kac formula,
\item  \cite{HolKarlKenLie10Splitting} for operator splitting methods for rough solutions, 
\item  \cite{GlowOshYin17Splitting, GloLeTa89Augmented, Glo03Finite, Ya71Method, RyuYin23large_scale, BauCom17convex, Mar82methods, Mar88splitting, Mar86splitting} for surveys, innumerable references, applications, and using operator splitting numerical methods for solving PDEs and variational nonlinear inequalities, (convex) optimization problems (including sparse and large scale problems), fixed point algorithms, problems for monotone operators etc.,
\item   \cite{HaiLu10geometric, BlaCa16concise, McLachQuis02splitting,SanzCal94numerical} for applications to geometric integration of ODEs. 
\end{itemize}
All these sources contain vast bibliographies. Additional references will be given later in the text.    

In what follows the aforementioned sources are used as starting points.

\section{Matrices and the Lie product formula}
 
We start with recalling the famous Lie (product) formula\footnote{also referred to as the Trotter (product) formula or the exponential product formula} for matrices. Let $M_n$ be the set of $n\times n$ matrices over $\mbC,$ $n\in\mbN.$ 
The function $y(t)=\e^{tA}y_0, t\in\mbR_+,$ is the solution to  
\[
\frac{dy(t)}{dt} = Ay(t)
\]
given $y(0)=y_0.$ 

\begin{theorem}[e.g. \cite{ReedSai72methods, Hall13quantum}]
\label{th:lie.product}
For $A_1,A_2\in M_n$
\begin{equation*}
\label{eq:1}
 \left\| \left(\e^{\frac{1}{n}A_1} \e^{\frac{1}{n}A_2}\right)^n-  \e^{A_1+A_2} \right\| \le \frac{C}{n},
\end{equation*}
where $\| \cdot\|$ is an arbitrary matrix norm and $C=C(A,B).$
\end{theorem}

\begin{remark}
The original source for Theorem \ref{th:lie.product} seems to be untrackable~\cite{CoFriedKaKe82eigenvalue}.
\end{remark}

One possible proof of the Lie product formula is based on the following telescopic identity: for $A,B\in M_m$ and any $n\in\mbN$
\begin{equation}
\label{eq:2}
A^n - B^n = \sum_{k=\ov{0,n-1}} A^k (A-B) B^{n-1-k},
\end{equation}
where $A^0=B^0=\Id$.

Let $S_n = \e^{\frac{1}{n}(A+B)}, T_n = \e^{\frac{1}{n}A}\e^{\frac{1}{n}B}.$ The identity \eqref{eq:2} implies 
\begin{align}
\label{eq:3}
\left\| \left(\e^{\frac{1}{n}A} \e^{\frac{1}{n}B}\right)^n - \e^{A+B} \right\| &\le n \left\| S_n - T_n \right\| \max\left\{\| T_n \|, S_n\| \|\right\}^{n-1} \\
&\le n\e^{\|A\| +\|B\|} \left\| S_n - T_n \right\|, \notag
\end{align}
where the local error is  
\begin{equation}
\label{eq:2.add.1}
\e^{\frac{1}{n}A} - \e^{\frac{1}{n}A_1}\e^{\frac{1}{n}A_2} = \frac{1}{2n^2} \left[ A_2, A_1\right] + O\left(\tfrac{1}{n^3}\right),
\end{equation}
so
\[
 \left\| S_n - T_n \right\|\le \frac{C}{n^2},
\]
which finishes the proof\footnote{We get the first order for the whole method since \eqref{eq:2.add.1} is a local error which gets propagated exactly $n$ times. This is a universal rule (cf.~\cite{Mil95numerical}).}.

\begin{remark}
\eqref{eq:2} is an arithmetic identity and thus holds for unbounded operators, too. Despite its triviality, \eqref{eq:2} is a standard and well-known tool for estimating the global error in the terms of local errors and  can be encountered on numerous occasions across the whole selection of works referenced in the present paper. 
\end{remark}

As we will see, the idea of the proof of Theorem~\ref{th:lie.product} extends into the abstract setting directly if a suitable form of uniform boundedness (''stability'') holds for $\max\left\{\| T_n \|, S_n\| \|\right\}^{n-1}$ in~\eqref{eq:3}. Otherwise, new methods should be proposed. 

\eqref{eq:2.add.1} implies that the method is at least of the second order if the matrices $A_1$ and $A_2$ commute, that is, if
\[
\left[A_1,A_2 \right] = A_1A_2- A_2A_1 = 0.
\]
However, if matrices commute, so do the corresponding exponentials and hence  the method is exact in fact ($C=0$ in the statement of Theorem~\ref{th:lie.product}). 

\begin{remark}
Alternative proofs of Theorem~\ref{th:lie.product} usually  use properties of the matrix logarithm additionally and cannot be easily extended into the general abstract setting.
\end{remark}

Now we can give the simplest example of  a operator splitting technique. 

\begin{example}
\label{ex:2.1}
Given $A_1, A_2 \in M_n$ consider ODEs
\[
\frac{dy(t)}{dt} = (A_1+A_2)y(t), \quad y(0)= y_0,
\]
and 
\begin{align*}
\frac{dy_k(t)}{dt} = A_k y_k(y), \quad k =1,2.
\end{align*}
Then, by Theorem~\ref{th:lie.product} for any $y_0$
\begin{equation*}
\label{eq:4}
\left\| \left(\e^{\frac{t}{n}A_1}\e^{\frac{t}{n}A_2}\right)^n y_0 - y(t) \right\| \le \frac{C}{n} \|y_0\|, \quad n\in\mbN.
\end{equation*}
\end{example}

\begin{remark}
Theorem~\ref{th:lie.product} also finds applications in the theory of matrix Lie groups and the corresponding homeomorphisms~\cite{Hall15Lie}, matrix trace inequalities~\cite{Petz94survey} etc.
\end{remark}

The next theorem was originally proposed for one two-dimensional difference scheme and can also be proved by direct calculations. 
\begin{theorem}[e.g.~\cite{Strang68construction}]
\label{th:strang.product}
For $A_1,A_2\in M_n$
\begin{equation*}
\label{eq:1.2}
 \left\| \left(\e^{\frac{1}{2n}A_1} \e^{\frac{1}{n}A_2} \e^{\frac{1}{2n}A_1}\right)^n-  \e^{A_1+A_2} \right\| \le \frac{C}{n^2},
\end{equation*}
where $\| \cdot\|$ is an arbitrary matrix norm and $C=C(A,B).$
\end{theorem}

\begin{remark}
There is a common agreement that G.~Strang and G.I.~Marchuk independently proposed to use the idea of Theorem~\ref{th:strang.product} to derive one very popular operator splitting method.
\end{remark}

Obviously, Theorem~\ref{th:strang.product} can be used to solve the ODE in Example~\ref{ex:2.1}, too:
\[
\left\| \left(\e^{\frac{t}{2n}A_1}\e^{\frac{t}{n}A_2}\e^{\frac{t}{2n}A_1}\right)^n y_0 - y(t) \right\| \le \frac{C}{n^2} \|y_0\|, \quad n\in\mbN.
\]
This is a simple installment of the Strang splitting scheme 
 which achieves the second order of accuracy with almost no increase in the number of necessary calculations since  three subsequent steps equal
\[
\Big(\e^{\frac{1}{2n}A_1} \e^{\frac{1}{n}A_2} \e^{\frac{1}{2n}A_1}\Big)^3 = \e^{\frac{1}{2n}A_1} \Big(\e^{\frac{1}{n}A_2}  \e^{\frac{1}{n}A_1}  \Big)^2 \e^{\frac{1}{n}A_2}  \e^{\frac{1}{2n}A_1}
\]
and so on.

The local error of the Strang splitting is~\cite{HundsVer07Numerical}
\begin{equation}
\label{eq:2.add.2}
 \e^{h(A_1+A_2)} - \e^{\frac{h}{2}A_1} \e^{hA_2} \e^{\frac{h}{2}A_1} = \frac{h^3}{24}\left( [A_2, [A_2, A_1]] + 2[A_1, [A_2, A_1]] \right) + O(h^4). 
\end{equation}

\section{Semigroups and the Chernoff product formula for abstract Cauchy problems} 
To venture forth, we need to recall standard material on semigroups and the abstract Cauchy problem (ACP), the Trotter-Kato theorems and the Chernoff product formula~\cite{Pazy12Semigroups, Goldstein85semigroups, EnNa00One-parameter, Da0One-parameter, ItoKa02evolution, Cher74product, EthKur09Markov}. We omit some details and technicalities. All integrals are  Bochner integrals.

Let $X$ be a Banach space. A family $(T_t)_{t\ge 0}$ of bounded linear operators on $X$ is a (one-parameter) strongly continuous  semigroup, or $C_0-$semigroup, if $T_0=\Id, T_{t+s}=T_t T_s, s,t \ge 0,$ and the mapping $t\mapsto T_t, t\ge 0,$ is strongly continuous. 

We denote the domain of a linear operator $A$ by $D(A).$

 A linear operator $A$ is the generator  of a $C_0-$semigroup $(T_t)_{t\ge 0}$ if
\begin{equation}
\label{eq:3.1}
Ax = \lim_{s\to 0+} \frac{1}{s}(T_s x - x)
\end{equation}
for all $x$ in 
\[
D(A) =\left\{ x\mid \mathrm{the\ limit\ in\ \eqref{eq:3.1}\ exists} \right\}.
\]

A generator is closed and densely defined. The pair $(A,D(A))$ defines the semigroup uniquely. 
For any $C_0-$semigroup $(T_t)_{t\ge 0}$ there exists $M \ge 1$ and $\omega \ge 0$ s.t.
\begin{equation*}
\label{eq:3.2}
\|T_t\| \le M \e^{\omega t}, \quad t\ge 0.
\end{equation*}
If $M=1$ and $\omega=0$ $(T_t)_{t\ge 0}$ is called a contraction semigroup. 

Though it is customary to write $(\e^{tA})_{t\in \mbR_+}$ for the semigroup generated by $A$ we do not follow this convention unless stated otherwise. 

\begin{example}
If $A$ is a bounded linear operator (e.g. $A\in M_n$) the family $\{\e^{tA}\mid t\ge 0\}$ is a $C_0-$semigroup with generator $A.$ 
\end{example}

\begin{example}
Any Feller $\mbR^m-$valued process $\{\xi(t)\mid t \ge 0\}$ automatically defines a contraction $C_0-$semigroup  
\[
T_t f(x) = \E_x f(\xi(t)), \quad f\in C_0(\mbR^m),
\]
where $C_0(\mbR^m)$ is the space of continuous function that vanish at infinity~\cite{EthKur09Markov, BottSchiWang13Levy}\footnote{The very definition of a Feller process may vary between authors~\cite{BottSchiWang13Levy}.}. In particular, Feller's one-dimensional diffusion on an interval $I$ has  generator 
\[
A = \frac{1}{2}a(x) \frac{d^2}{dx^2} + b(x) \frac{d}{dx} + c(x),
\]
while the precise description of $D(A)$ incorporates boundary conditions (and thus the behavior of the process at the boundary) and $c$ corresponds to killing. Note that it is not possible in general to give a full description of  $D(A)$ for a Feller process 
and it is often sufficient to consider a core $L$ instead, that is, a set $L\subset D(A)$ such that the closure of $A\vert_{L}$ equals $A$~\cite{EthKur09Markov}. 
\end{example}

\begin{remark}
However, the heat semigroup is not a $C_0-$semigroup on $C_b(\mbR^m).$ Thus one should carefully choose a space that a semigroup acts on.
\end{remark}

\begin{example}
A Wiener process analogously  defines a $C_0-$semigroup on $L_2(\mbR^m)$ with generator $\frac{1}{2}\Delta$ and $D(A) = W^{2,2}(\mbR^m)$~\cite{BoKryRo15Fokker}. However, the possibility to extend a Feller semigroup onto some $L_p$ is linked to the properties of the adjoint operator and thus to properties of Fokker-Plank-Kolmogorov equations (for discussions and conditions see e.g.~\cite{Stroock2008Partial, Cas85generators, BottSchiWang13Levy} and \cite[Section VI.4]{EnNa00One-parameter}).
\end{example}

Assume that $A$ is a densely defined closed linear operator with a non-empty resolvent set on a Banach space $X.$ We  consider the  autonomous inhomogeneous abstract Cauchy problem 
\begin{align}
\label{eq:3.3}
\frac{dy(t)}{dt} &= Ay(t) + f(t), \quad t \in (0;T) \notag \\
y(0) &= y_0. 
\end{align}

Given $y_0\in D(A)$ a classical  solution to \eqref{eq:3.3} is a $X-$valued function $u$ s.t. $y(t)\in D(A), t\in [0;T],$ $y$ is continuous on $[0;T]$ and continuously  differentiable (in the terms of the norm on $X$) on $(0;T),$ and \eqref{eq:3.3} holds.   A homogeneous ACP is well-posed if it has the unique classical solution that depends (uniformly) continuously on the initial data. 

The following result gives the well-known relation ACPs and $C_0-$semigroups.

\begin{theorem}
\label{th:3.1}
\begin{enumerate}
\item 
The homogeneous ACP \eqref{eq:3.3} is well-posed if and only if $A$ generates a $C_0-$semigroup. 
\item If $A$ generates a $C_0-$semigroup $(T_t)_{t\ge 0},$ $y_0\in D(A),$ the function $f\in L_1((0;T))\cap C((0;T))$ takes values in $D(A)$ and $Af \in L_1((0;T))$  then the function
 \begin{equation}
\label{eq:3.4}
y(t) = T_t y_0 + \int_0^t T_{t-s} f(s) ds, \quad t \in[0;T],
\end{equation}
 is the unique classical solution of \eqref{eq:3.3}.
\end{enumerate}
\end{theorem}

There are other types of solution. In particular, the mild solution \eqref{eq:3.4} is often used when the functions $f$ and $y$ fail to satisfy regularity assumptions stated above. 

The non-autonomous ACP 
\begin{align}
\label{eq:3.a.1}
\frac{dy(t)}{dt} &= A(t)y(t) + f(t), \quad t \in (s;T) \notag \\
y(s) &= y_0,
\end{align}
is  called an evolution problem. A classical solution is defined similarly to the autonomous case but no simple analog of Theorem~\ref{th:3.1} exists (see \cite{Pazy12Semigroups, Goldstein85semigroups} for the theory and \cite{NaNi02well-posedness, Schnau02wellposedness}  and \cite[Remark 1.12]{BatCsoFarNi11operator} for a survey and references).  Moreover, \eqref{eq:3.a.1} requires the compatibility of $\{D(A(t))\mid t \in [0;T]\}$ if the domains in question are time-dependent, a standard condition being $\ov{\cap_{t\in[0;T]} D(A(t))}= X.$ The equation \eqref{eq:3.a.1} is then solved for $x\in \cap_{t\in[0;T]} D(A(t)).$ 

However, the well-posedness of a non-autonomous homogeneous ACP implies the existence of an evolution family of bounded operators (also called propagators) $(U_{s,t})_{0\le s \le t \le T}$ s.t. $U_{t,t}=\Id, U_{s,t}U_{t,r} =U_{s,r},  s \le t\le r,$ and the mapping 
$(s,t) \mapsto U_{s,t}$  is strongly continuous. For any fixed $s\in [0;T]$ the solution for 
\begin{align*}
\frac{dy(s;t)}{dt} &= A(t) y(s;t), \quad t \in (s;T), \\
y(s;s) &= y_0
\end{align*}
is then given via
\[
y(s;t) = U_{s,t} y_0, \quad t\in[s;T],
\]
and a mild solution to the analog of \eqref{eq:3.a.1} with $f\in L_1((0;T))$ is given via
\[
y(s;t) = U_{s,t} y_0 + \int_{s}^t U_{t,r} f(r) dr, \quad t\in [s;T],
\]
for any $s\in[0;T].$ 

\begin{remark}
In general, non-autonomous inhomogeneous ACPs, unavoidable in applications, may behave in an unexpected way and one should be careful to draw conclusions about the properties of the corresponding propagators~\cite{NaNi02well-posedness}. 
\end{remark}

Until further notice, we consider only homogeneous ACPs.

The following theorem is a version of the celebrated Trotter-Kato theorem\footnote{occasionally called the  Trotter-Neveu-Kato theorem (e.g~\cite{Zag22notes, CaZag01operator, Bob07limitations})}  whose general form provides a basis for numerous approximation schemes in the theory of semigroups, including  Yosida approximations\footnote{For the exposition of SPDEs in the terms of Yosida approximations see~\cite{Go16Yosida}.}. One can consult \cite{GoTo14convergence, GoKoTo19general, ItoKa02evolution}  for a survey of the general approximation theory of $C_0-$ semigroups, some results about the rate of the convergence (in the operator norm, too) and further references.  A number of references concerning the Chernoff and Trotter-Kato product formulas will be given later.

\begin{remark}
\cite{Pfei85approximation,Pfei84probabilistic,Pfei86some,Pfei82general} use probabilistic  methods to develop a unified exposition of many approximation formulas (see also \cite{BentPau04optimal}). 
\end{remark}

Roughly speaking, Trotter-Kato theorems connect convergence of semigroups, their resolvents and generators. 

\begin{theorem}
\label{th:a.1}
Let $(T_{nt})_{t\ge 0}, n\in\mbN,$ be $C_0-$semigroups on a Banach space $X$ with generators $(A_n,D(A_n)), n\in\mbN,$ s.t. 
\begin{equation}
\label{eq:3.stability}
\|T_{nt}\| \le M \e^{\omega t}, \quad t\ge0, n\in\mbN,
\end{equation}
for some $M\ge 1, \omega\in\mbR.$ For fixed $\lambda > \omega$ and the following assertions 
\begin{enumerate}
\item There exists a densely defined operator $(A,D(A))$ s.t. $A_n\to A, n\to\infty,$ strongly on a core of $A$ and the range of $\lambda\Id-A$ is dense in $X;$
\item There exists a bounded linear operator $R$ with dense range s.t. $(A_n - \lambda \Id)^{-1}\to R, n\to\infty,$ strongly;
\item The semigroups  $(T_{nt})_{t\ge 0}, n\in\mbN,$ converge strongly and uniformly on bounded sets to a semigroup with generator $B$ s.t. $R=(B-\lambda \Id)^{-1};$ 
\end{enumerate} 
it holds that (1) $\Rightarrow$ (2) and (2) $\Leftrightarrow$ (3).
\end{theorem}

\begin{remark}
A probabilistic formulation of Theorem~\ref{th:a.1} for Feller processes can be found in~\cite{kallenberg}. Trotter-Kato theorems are behind some results about approximations of Feller processes with Markov chains, in particular~\cite{kallenberg}.
\end{remark}

\begin{remark}
Counterexamples for the Trotter-Kato theorem, in particular of probabilistic origin, can be found in~\cite{Bob07limitations};  variational and nonlinear versions of the Trotter-Kato theorem, in~\cite{ItoKa02evolution}; generalizations for norm operator topology, in \cite{Zag22notes, CaZag01operator}.
\end{remark}

The next  theorem is called the Chernoff product formula and is another fundamental tool in the approximation theory for $C_0-$semigroups.

\begin{theorem}
\label{th:3.4}
Let $V(t), t\ge 0,$ be a family of bounded linear operators s.t. $V(0)=\Id$ and 
\begin{equation}
\label{eq:3.stability.ch}
\|V(t)^n\| \le M\e^{\omega n t}, \quad t \ge 0, n\in\mbN, 
\end{equation}
for some $M\in\mbR_+, \omega\in\mbR,$ and let the limit
\[
Ax = \lim_{t\to 0} \frac{V(t)x - x}{t}
\]
exist for all $x\in D$ where $D$ and $(\lambda\Id - A)D$ are dense subspaces of $X$ for some $\lambda > \omega.$ Then the closure of $A$ generates a $C_0-$semigroup $(T_t)_{t\ge 0}$ and for any sequence of positive integers $(k_n)_{n\in\mbN}$ and any positive sequence $(t_n)_{n\in\mbN}$ satisfying $k_n t_n \to t, t_n \to 0, n\to \infty,$
\begin{equation}
\label{eq:3.a.2}
T_t x = \lim_{n\to\infty} (V_{t_n})^{k_n} x
\end{equation}
for any $x\in D.$ If $t_n k_n =t, n\in\mbN,$ the convergence is uniform on bounded intervals.  
\end{theorem}

The formula  \eqref{eq:3.a.2} provides an approximation of the semigroup $(T_t)_{t\ge 0}$ in the terms of the family $V.$ 
The exposition of the method of obtaining  such approximations is given in a survey~\cite{But20Method}; further references and applications to stochastic processes including those on manifolds and in domains with Dirichlet and Robin boundary conditions  can be found in~\cite{But20Method, But18Chernoff, But18Chernoff_for_killed, MaMoReSmo23Chernoff, ButGroSmo10Lagrangian, Nitt09approximations, Sha11Chernoffs, Ne09ChernoffArx, SmoWeizWi07Chernoffs, MaMoReSmo23Chernoff, Ob06representation, SmoToTru02Hamiltonian, But08Feynman, ButSchiSmo11Hamiltonian, ButGroSmo10Lagrangian, But18Chernoff, SmoWeizWi00Brownian}. Other references given later in the text for the Trotter-Kato theorem in the context of the approximation theory are relevant, too.

The following additional references illustrate a possibility to strengthen the conclusion of the Chernoff product formula:
\cite{Zag22notes,Zag17comments, Za20notes_chernoff_formula, CaZag01operator,GalRe22rate, GalRe21upper_arx, Prud20speedArx} establish  convergence in the operator norm  an/or  obtain estimates of the rate of the convergence  for self-adjoint operators and quasi-sectorial contraction semigroups etc.; \cite{Pau04operator-norm, Zag17comments, Zag22notes, CaZag01operator}  obtain estimates by using probabilistic arguments. The convergence can be arbitrary slow and may not hold in a stronger topology. 

\begin{example}
Theorem~\ref{th:3.4} yields the exponential formula
\begin{equation}
\label{eq:3.euler.b}
T_t = \lim_{n\to \infty}\Big(\Id - \frac{t}{n}A\Big)^{-n}.
\end{equation}
\end{example}

\begin{example}
If $A$ is a bounded operator Theorem~\ref{th:3.4}  can be used to deduce
\begin{equation}
\label{eq:3.euler.f}
\e^{tA} =\lim_{n\to \infty} \left(\Id + \frac{t}{n}A \right)^n.
\end{equation}
if the series representation of $\e^{tA}$ is taken as a definition of the exponential. 
\end{example}

\begin{example}\label{ex:1}(\cite{But20Method}; cf. \cite{Lejay18Girsanov}) 
Assume that real-valued functions $a, \sigma$ are Lipschitz continuous and bounded, $\inf_{u\in\mbR} |\sigma(u)| > 0$ and consider an SDE
\[
dx(t) = a(x(t))dt  + \sigma(x(t))dw(t),\quad t\in [0;T],
\] 
where $w$ is a standard Wiener process. The Euler-Maruyama approximations  are 
\begin{align*}
y_{n,k+1} &=  y_{n,k} + a(y_{n,k})\frac{T}{n} + \sigma(y_{n,k}) \Big(w\Big(\frac{T(k+1)}{n}\Big) - w\Big(\frac{Tk}{n}\Big) \Big), \\
& k = \ov{0,n-1}, n\in\mbN.
\end{align*}
For any $f\in Lip(\mbR)$ 
\[
 \E f(y_{n,n})  \to  \E f(x(T)) , \quad n\to \infty.
\]
On the other side, one can  show that the family $(V(t))_{t\ge 0}$ defined via
\[
V(t)f(x) = \frac{1}{(2\pi t)^{1/2}\sigma(x)} \int_{\mbR} \e^{-\frac{(y-x-a(x)t)^2}{2t\sigma(x)^2}}  f(y) dy, \quad t \ge 0,
\]   
satisfies the conditions of Theorem~\ref{th:3.4} implying 
\begin{equation}
\label{eq:3.a.5}
\E f(x(T)) = \lim_{n\to\infty} V\Big(\frac{T}{n}\Big)^n f(x)
\end{equation}
for sufficiently smooth compactly supported $f.$ See also \cite{SmoWeizWi07Chernoffs, MaMoReSmo23Chernoff, SmoWeizWi00Brownian} for other applications of the Chernoff product formula to pathwise approximations of random processes. Note that this result does not recover the first order of the weak Euler-Maruyama scheme.
\end{example}

\begin{example}\label{ex:2}(\cite{ButGroSmo10Lagrangian}) Let $g_{\mu}$ be the one-dimensional Gaussian density with $0$ mean and variance $\mu.$ \eqref{eq:3.a.5} can be rewritten as
\begin{align}
\label{eq:3.a.6}
\E_{x_0} f(x(T)) &= \lim_{n\to\infty} \int_{\mbR^n} f(x_n) \exp\Bigl\{ \sum_{k=\ov{1,n}} \frac{a(x_{k-1}) (x_k -x_{k-1}) -\frac{t }{2n}a(x_{k-1})^2   }{\sigma(x_{k-1})^2} \Bigr\} \notag \\
&\qquad\qquad \times \prod_{k=\ov{1,n}} g_{\frac{T\sigma(x_{k-1})^2}{n}}(x_k-x_{k-1}) dx_1\ldots dx_n.
\end{align}
Such a representation of the original semigroup as the limit of iterated integrals w.r.t. some finite-dimensional projections  of a (pseudo)measure on a phase space (the Wiener measure on $C([0;T])$ here) is called the Feynman formula~\cite{SmoToTru02Hamiltonian}\footnote{The definition of the Feynman formula may slightly vary (c.f.~\cite{ButSchiSmo11Hamiltonian}).} and presents an alternative approach to the Feynman-Kac formula.  
See \cite{But20Method, But18Chernoff_for_killed, SmoWeizWi07Chernoffs, MaMoReSmo23Chernoff, Ob06representation, SmoToTru02Hamiltonian, But08Feynman, ButSchiSmo11Hamiltonian, ButGroSmo10Lagrangian, But18Chernoff, SmoWeizWi00Brownian}  for other results on Feynman formulas.
\end{example}

\begin{example}\label{ex:a.1}\cite{But20Method}
Setting 
\[
dy(y) =\sigma(y(t))dw(t), \quad t \in[0;T],
\]
and passing to the limit in~\eqref{eq:3.a.5} and~\eqref{eq:3.a.6} gives a particular case of the Girsanov theorem~(e.g.~\cite{LipShi77Statistics}) for $x$ and $y:$
\[
\E_{x_0} f(x(T)) = \E_{x_0} f(y(T)) \exp\Big\{\int_0^T \frac{a(y(t))}{\sigma(y(t))}dw(t) - \frac{1}{2}\int_0^T \frac{a(y(t))^2}{\sigma(y(t))^2} dt\Big\}.
\] 
\end{example}

\begin{remark}
The behavior of a solution of an PDE and its differentiability at $t=0$ (and thus the behavior of the corresponding semigroup) can be a subtle moment in applications, particularly when a probabilistic interpretation is used. A typical example is the function $u(x,t)=\E_x f(\xi(\min\{t,\tau\})),$ where $\tau$ is the moment a Feller process $\xi$ hits the  boundary of a domain $D.$ If $f=\1_{\pt D}$ then $u= \Prob(\tau <t)$ is discontinuous at $t=0$ though $u\in C^{2,1}(D^{int}\times (0;\infty))$ usually.  As a result, the weak variational formulation in the terms of Gelfand triples~\cite{PaWin08First} or the weak distributional formulation~\cite{LeRaRey22probabilistic, TriZab11Pfaffian} can be suggested as an alternative (see also~\cite{FeePop15stochastic, Wang17viscosity}).    
\end{remark}  

\begin{remark}
One particular example of the versatility of the Chernoff product formula is that it can be used to prove the central limit theorem~\cite{Goldstein85semigroups} .
\end{remark}

\section{The Trotter-Kato formula}

Now we state one of core results of the theory of abstract operator splitting, the Trotter-Kato  formula.

\begin{theorem}
\label{th:3.t-k}
Let $(T_t)_{t\ge 0}$ and $(S_t)_{t\ge 0}$ be $C_0-$semigroups on a Banach space $X$ with generators $(A,D(A))$ and $(B,D(B)),$ respectively, s.t.
\begin{equation}
\label{eq:3.stability.tk}
\|T_t^n S_t^n\| \le M\e^{\omega nt}, \quad n\in\mbN, t \ge 0,
\end{equation}
for some $M\ge 1$ and $\omega\in\mbR.$ Define $C=A+B$ on $D=D(A) \cap D(B).$ If $D$ and $(\lambda\Id- A -B)D$ for some $\lambda > \omega$ are dense in $X,$ then the closure  of $C$ generates a $C_0-$semigroup $(U_t)_{t\ge 0}$  and 
\[
U_tx = \lim_{n\to\infty} \left(T_{t/n}S_{t/n}\right)^n x
\]
uniformly on compact intervals in $t.$
\end{theorem}

Theorem~\ref{th:3.t-k} can be formulated for a finite number of semigroups (e.g.~\cite{Pazy12Semigroups}).

\begin{remark}
This result belongs to H.E.~Trotter and T.~Kato. Yu.L.~Daletskii also obtained similar product formulae at the same time (see \cite{DaFo91measures} for references).
\end{remark}

\begin{remark}
The choice of the topology matters: the convergence may not hold in a stronger topology~\cite{NeidSteZa18remarks}.
\end{remark}

\begin{remark}
The stability assumption~\eqref{eq:3.stability.tk} cannot be dropped~\cite{KuhnWa00Lie}.
\end{remark}

\begin{example}
\label{ex:4.4}
Assume that $f$ is a complex-valued function and $\int_{\mbR^d} |f(x)|^2 dx=1.$  The Shr\"{o}dinger operator 
\[
H = \frac{1}{2}\Delta + V
\]
describes the motion of a (spinless) quantum particle under the action of the real-valued potential $V$ as follows. The wave function  $u$ is the solution of 
\begin{align}
\label{eq:3.9}
\frac{\p u}{ \p t}  &= iHu, \notag \\
u\vert_{t=0} &= f,
\end{align}
and the probability density at time $t$ of the position of the particle is $|u(x,t)|^2.$ \eqref{eq:3.9} can be interpreted as an ACP 
\[
\frac{du(t)}{dt} = \Big(\frac{i}{2}\Delta + i V\Big) u(t)
\]
in some Hilbert space.  The operator $H$ is essentially self-adjoint (that is, has a self-adjoint extension) under some rather mild assumptions on $V,$ so the semigroup generated by $iH$ is an unitary  $C_0-$semigroup, in particular\footnote{We skip all technical aspects associated with the existence of multiple extensions etc., and we do not mention symmetric and skew-adjoint operators even though treating these issues is an integral part of quantum mechanics.}. Since the  semigroups for $\frac{i}{2}\Delta$ and  $iV$ are
\begin{align*}
S_t f(x) &= \frac{1}{\left(2 i \pi t\right)^{d/2}} \int_{\mbR^d} \e^{\frac{i \|x-y\|^2}{2t}}  f(y) dy, \\
T_t f(x) &= \e^{it V(x)} f(x), \quad t \ge 0,
\end{align*}
respectively, the Trotter-Kato theorem implies 
\begin{align*}
u(x_0,t) &= \lim_{n\to\infty}\left(S_{t/n} T_{t/n}\right)^n f(x_0) \\
&=  \lim_{n\to\infty} \frac{1}{\left(2 i \pi t\right)^{nd/2}} \int_{\mbR^{nd}} \exp\Big\{ \frac{in}{2t} \sum_{k=\ov{0,n-1}} \|x_k - x_{k+1}\|^2  + \frac{it}{n} \sum_{k=\ov{0,n-1}} V(x_{k+1})   \Big\} \\ 
&\qquad\qquad\qquad\qquad\qquad \times f(x_n) dx_1 \ldots dx_n,
\end{align*}
which is an polygonal approximation for the Feynman integral over all histories (trajectories in the space $C_x([0;t])$ of continuous trajectories starting at $x$)~\cite{ReedSai72methods, ReedSi75methods_2, Goldstein85semigroups, JohnLa00Feynman}
\begin{align}
\label{eq:3.10}
C \int_{C_x([0;T])} \e^{iS(\omega)} f(\omega(s)) \prod_{0\le s \le t} d\omega_s,
\end{align}
where 
\[
S(\omega) = \int_0^t \Big( \frac{1}{2} \Big\| \frac{d\omega(s)}{ ds} \Big\|^2 + V(\omega(s))\Big) ds
\]
is the action integral. The expression~\eqref{eq:3.10} is purely formal and the Trotter-Kato formula is a standard way to give this representation a rigorous meaning. 
\end{example}

\begin{example}(\cite{ReedSai72methods, Si05functional, JohnLa00Feynman, LoHiBetz20FeynmanKac, Hall13quantum, ReedSi75methods_2})
It is customary for textbooks on quantum mechanics to derive the  Feynman-Kac formula for Shr\"{o}dinger operators as a corollary of the Trotter-Kato formula. To understand the idea formally recall that the solution to the following autonomous homogeneous PDE (both assumptions are used in what follows)
\begin{align*}
\frac{\pt u(x,t)}{\pt t} &= \frac{1}{2}\Delta u(x,t) + V(x) u(x,t), \quad t \ge 0, x \in\mbR^d, \\
u(x,0) &= f(x),\quad  x \in \mbR^d, 
\end{align*}
is given via
\[
u(x,t) = \E f(x+w(t)) \e^{\int_{0}^t V(x+w(s))ds}, 
\]
where $w$ is a standard Wiener process. 
On the other hand, the solution to 
\[
\frac{\pt q(x,t)}{\pt t} = V(x) q(x,t)
\]
is 
\[
q(x,t) = S_t q(x,0) = \e^{V(x)t} q(x, 0),
\]
so combining this semigroup with the heat semigroup $(T_t)_{t\ge 0}$ in accordance with the Trotter-Kato theorem yields the approximations 
\begin{align*}
\Big( S_{1/n} \circ T_{1/n}\Big)^n f(x) &= \E f(x + w(t)) \exp\Bigl\{ \frac{t}{n} \sum_{k=\ov{0,n-1}} V\Bigl(x + w(t) - w\Bigl(\frac{kt}{n}\Bigr)\Bigr)\Bigr\} \\ 
&= \E f(x+w(t)) \exp\Big\{ \frac{t}{n} \sum_{k=\ov{1,n}} V\Big(x+w\Big(\frac{kt}{n}\Big)\Big)\Bigr\}, \\
&\qquad\qquad n \in\mbN,
\end{align*}
of the solution $u$ (similarly for $( T_{1/n} \circ S_{1/n})^n)$). 
\end{example}

\begin{remark} 
It should be noted that precise formulations of  physical and probabilistic versions of the Feynman-Kac formula may be dictated by different points of focus:  physicists are often more concerned with results for the Laplacian and singular or irregular (in different senses) potentials (e.g. those belonging to Kato classes or $L_\infty(\mbR^d)$)\footnote{On this path, various approximations of the original problem with ACPs with ''better'' coefficients are needed along with the corresponding limit theorems.}. This remark applies primarily to textbooks.
\end{remark}

\begin{example}(\cite{Lejay18Girsanov}; cf. Example~\ref{ex:a.1}) It is possible to derive the Girsanov theorem using the Trotter-Kato formula (under some assumptions on the regularity of drift)\footnote{\cite{Lejay18Girsanov} also uses the Euler-Maruyama scheme to prove the Girsanov theorem.}. Following the original source, we explain some intuition behind this. Let  $(T_t)_{t\ge 0}$ be the heat semigroup. The semigroup $(S_t)_{t\ge 0}$
\[
\frac{\pt S_t f(x)}{\pt t} = a(x) \cdot \nabla S_tf(x) 
\]
  is alternatively given for small $t$ as
\begin{align*}
S_tf(x) &= f(\xi_x(t)) = f(x) + \nabla f(x)\cdot (\xi_x(t) -x)+ o(t) \\
&=  f(x) + t \nabla f(x)\cdot  a(x) + o(t),
\end{align*}
where $\xi_x(t) = x + \int_0^t a(\xi_x(s)) ds.$ It can thus be shown that 
\begin{align*}
\left( T_t \circ S_t\right) f(x) &= \E \left( f(x+w(t)) + ta(x)\cdot \nabla f(x + w(t)) \right)\\ 
&\qquad  + t \E \sum_{k=\ov{1,m}} \frac{\pt f(x + w(t))}{\pt x_k}\frac{\pt a(x)}{\pt x_k} w_k(t) + o(t) \\ 
&= \E f(x + w(t)) (1 + a(x)\cdot w(t)) + o(t) \\
&= \E f(x + w(t)) \e^{a(x) \cdot w(t) -\frac{t}{2} \|a(x)\|^2} + o(t),
\end{align*}
for a Wiener process $w$ since\footnote{This is a calculation from the Malliavin calculus. See also Section~11.}
\[
\E t \nabla f(x + w(t)) = \E  f(x + w(t)) w(t). 
\]
\end{example}

\begin{example}
The Trotter-Kato formula for self-adjoint operators~\cite{ReedSai72methods, Hall13quantum} is usually written\footnote{Here   $(\e^{tC})_{t\ge 0}$ denote the semigroup with generator $C$.} as
\[
\lim_{n\to\infty} \left(\e^{\frac{i}{n}A_1} \e^{\frac{i}{n}A_2} \right)^n = \e^{iA}
\]
where  $A_1, A_2$ are self-adjoint operators and  $A$ is the self-adjoint extension of $(A_1+A_2, D(A_1)\cap D(A_2))$. If $A_1$ and $A_2$ are  bounded below additionally: 
\[
\inf_{x\in D(A_k)} \frac{(x,A_k x)}{\|x\|^2} > -\infty, \quad k=1,2,
\]
 we also have
\[
\lim_{n\to\infty} \left(\e^{-\frac{1}{n}A_1} \e^{-\frac{1}{n}A_2} \right)^n = \e^{-A}.
\] 
The importance of these versions  for physics follows from the Stone theorem:  on Hilbert spaces,   every unitary $C_0-$group  has a generator of the form $iA$ for some self-adjoint $A$ and every symmetric $C_0-$semigroup has a generator of the form $-A$ for some bounded below self-adjoint $A,$ and this is exactly the situation of quantum mechanics.  The Trotter-Kato formula admits special proofs independent of the Trotter-Kato theorems in this case. In particular, if $A_1+A_2$ is densely defined and self-adjoint on $D(A_1)\cap D(A_2)$ the proof is rather direct. For further references and extensions see  \cite{LoHiBetz20FeynmanKac, JohnLa00Feynman, Si05functional}, in particular. However, the rate of the convergence was found only recently and under additional assumptions (see e.g. Example~\ref{ex:0.1} and related references).
\end{example}

\begin{example}\label{ex:2.non} To formulate a formal non-autonomous version of the Trotter-Kato formula, consider operators  $A_0, A_1, A_2$ such that the corresponding non-autonomous ACPs are well-posed with propagators $(T^{(k)}_{s,t})_{0\le s \le t \le T}, k =0,1,2,$ respectively. Additionally, assume that $D(A_0(t)),$ $D(A_1(t)),$ $D(A_2(t)),$  do not depend on $t,$ and $A_0(t)=A_1(t)+A_2(t)$ on $D(A_0(0)), t\in[0;T].$ Then one expects for any $s,t\in [0;T], s \le t,$ 
\begin{align*}
\lim_{n\to\infty} T^{(1)}_{s + \frac{(n-1)t}{n}, s+ t}\circ T^{(2)}_{s + \frac{(n-1)t}{n}, s+ t}\circ \ldots \circ T^{(1)}_{s, s+ \frac{t}{n}} \circ T^{(2)}_{s, s+ \frac{t}{n}}   = T^0_{s,t}.
\end{align*} 
Note that the condition of the domains being constant is a rather limiting one and excludes time-dependent boundary conditions in general as BC is usually incorporated into the very description of a domain. 
\end{example}


\begin{example}(\cite{EnNa00One-parameter}; also~\cite{BatCsoFarNi11operator})
The Trotter-Kato formula can be used to derive an approximation of a non-autonomous ACP by autonomous ACPs with  constant coefficients (cf.~Example~\ref{ex:2.non}). For that, assume that $(U_{s,t})_{0\le s \le t \le T}$ is the propagator for some well-posed non-autonomous ACP
and define semigroups $(T^{(s,n, k)}_{l})_{s \le l\le T}, k =\ov{1,n}, n\in\mbN, s \in[0;T],$ for the equations
\begin{align*}
\frac{dy_{s,n,k}(l)}{dl} &= A_{s,n,k} y_{s,n,k}(l), \quad l \in[0;T], \\
y_{s,n,k}(0) &= y_0, \\
A_{s,n,k} &=  A\Big(s + \frac{kt}{n}\Big), 
\end{align*} 
so that $y_{s,n,k}(l) = T^{(s,n,k)}_{l}y_0, l\in[0;T].$ Then 
\[
\lim_{n\to\infty}  T^{(s,n, n)}_{\frac{t}{n}} \circ T^{(s,n, n-1)}_{\frac{t}{n}} \circ \ldots \circ T^{(s,n,1)}_{\frac{t}{n}}   = U_{s,s+t} 
\]
strongly on $\cap_{t\in[0;T]} D(A_t)$ uniformly in $(s;t)$ on compact sets. 
\end{example}

\begin{remark}
In Example~\ref{ex:4.4},   $D(\Delta) \cap D(V\cdot)$ may be too small or just $\{0\}$ even if the sum $\frac{1}{2}\Delta +V$ can be defined via quadratic forms~\cite[Chapter 10]{JohnLa00Feynman}\cite[Chapter 9]{Hall13quantum}. In fact, one can use the Trotter-Kato formula to define the generalized sum of two operators with incompatible domains~\cite[Remark 8.17]{Goldstein85semigroups}.
\end{remark}

\begin{remark}
The original formulation of the Trotter-Kato theorem for PDEs requires boundary conditions to be time-independent.  In general, adding non-autonomous boundary and initial conditions can lead to  significant technical issues.
\end{remark}

Up to this moment, almost no attention has been given to the question of the rate of convergence even though it is an extremely important question for applications. Another feature  of the formulation of Theorem~\ref{th:3.t-k} is the usage of the strong topology. We know that either question may not have stronger results due to the references given for the Chernoff product formula: the rate of convergence can be arbitrary slow and convergence in stronger topologies may not hold. To understand the reasons behind these facts and thus to emphasize some limitations of the standard version of the Trotter-Kato formula we need to have a short discussion of the proof of Theorem~\ref{th:3.4} since Theorem \ref{th:3.t-k} is a direct corollary of Theorem~\ref{th:3.4}. 

Subsequent reasoning is well known. 

\begin{remark}
Note that the original proof of Theorem~\ref{th:lie.product} falls apart as soon as operators involved are unbounded (e.g. they contain differentiation). 
\end{remark}

The condition for the sets $D$ and $(\lambda\Id - A)$ to be dense in $X$ implies (in a nontrivial way) that a densely defined closed $A$ exists and is  indeed a generator of a $C_0-$semigroup so this condition can be seen~\cite{EnNa00One-parameter} as an extension of the Hille-Yosida  theorem\footnote{the implications $(1) \Rightarrow (2), (3)$ in Theorem~\ref{th:a.1}}. In fact, some formulations of the Trotter-Kato formula just require  $A+B$ to be a generator of a $C_0-$semigroup. The same applies to the Chernoff product formula (e.g.~\cite{Da0One-parameter}). In practical terms, this means that one has to be careful with domains when trying to replicate the original proof. 

The rest of the proof of Theorem~\ref{th:3.4} consists of two separate claims. The first one relies on the following observation: if $A$ is a linear bounded operator satisfying 
\[
\sup_{n\in\mbN} \|A^n\| \le M
\]
and $\xi$ is a Poisson random variable with mean $n\in\mbN$ then for any $x\in X$ 
\begin{align}
\label{eq:3.sqrt.n}
\Big\| \Big(A^n - \e^{n(A-\Id)}\Big)x\Big\| &\le \Var(\xi)^{1/2}   \|(A-\Id)x\| = n^{1/2} \|(A-\Id)x\|,
\end{align} 
and therefore  for bounded operators
\[
A_n = \frac{n}{t}\Big(V\Big(\frac{t}{n}\Big)- \Id\Big), \quad n\in\mbN,
\]
one can prove
\begin{align*}
\label{eq:3.5}
\Big\|\Big((T_t -V\Big(\frac{t}{n}\Big)\Big)x \Big\| &\le \Big\|\Big(T_t -\e^{tA_n}\Big)x \Big\| + \Big\|\Big(\e^{tA_n} -V\Big(\frac{t}{n}\Big)\Big)x \Big\| \notag \\
&\le  \Big\|\Big(T_t -\e^{tA_n}\Big)x \Big\| + \frac{t}{n^{1/2}} \Big\| \frac{V(\frac{t}{n}) -\Id}{t/n} x\Big\|, 
\end{align*}
where  
\[
\lim_{n\to\infty}\Big\| \frac{V(\frac{t}{n}) -\Id}{t/n} x\Big\|= \|Ax\|
\]
 by assumption. Thus it is left to prove
\begin{equation}
\label{eq:3.6}
 \left\|\left(T_t -\e^{tA_n}\right)x \right\| \to 0, \quad n\to \infty.
\end{equation}
This is the second part of the proof. Since $A_n \to A, n\to\infty,$  the Trotter-Kato theorem implies \eqref{eq:3.6}. But there are no possibilities to either have a bound for this expression or to replace strong convergence with convergence in the operator norm in general because the proof of the Trotter-Kato theorem does not produce such results internally. 
As a result, the standard Chernoff product formula and Trotter-Kato formula are formulated in the terms of strong convergence and with no bounds. 

It may be possible to refine~\eqref{eq:3.sqrt.n} and/or strengthen the convergence in the Trotter-Kato theorem for special classes of generators (e.g. self-adjoint or maximal accretive ones). Then one can combine  improved estimates in~\eqref{eq:3.sqrt.n} (in particular, those obtained by using properties of the Poisson distribution)  and consistency estimates for the Trotter-Kato theorem to yield  the conclusion about the convergence in the operator norm and to provide estimates for its speed~\cite{CaZag01operator, Zag22notes, Pau04operator-norm, Zag17comments, Za20notes_chernoff_formula} (see also~\cite[Chapter 3, Lemma 3.5]{Pazy12Semigroups}). This may be compared with \cite{GoTo14convergence, GoKoTo19general} where estimates of the speed of convergence often include an original element  $x.$ See also~\cite{ItoKa02evolution, ItoKa98Trotter} for other results of this type.  Some other references concerning bounds for the rate of convergence will be added later. 

The decomposition of the proof and conditions of Theorems~\ref{th:3.4} and \ref{th:3.t-k} is an illustration of one overarching principle of numerical analysis: stability and consistency imply convergence. This fundamental result is known as the Lax-Richtmyer theorem (the Lax equivalence theorem) in the context of finite difference approximations of PDEs:

\begin{theorem}
A consistent finite difference scheme for a well-posed linear initial value problem is convergent if and only if it is stable.
\end{theorem} 

This principle extends to other settings and situations as soon as an approximation scheme arises (see \cite{ItoKa02evolution} for a discussion in the setting of evolution equations and \cite[Chapter 3]{GlowOshYin17Splitting} for a discussion in the case of the Trotter-Kato formula). 

Conditions \eqref{eq:3.stability}, \eqref{eq:3.stability.ch} and \eqref{eq:3.stability.tk} give stability of a scheme while  the convergence  of generators  or the behavior of the family $(V(t))_{t\ge 0}$ at $0$ are consistency conditions in our case etc. However, one should not forget that, contrary to the case of the Lax-Richtmyer theorem, the principle ''stability and consistency imply convergence'' is no longer a rigorous statement and does not also reflect the full nature of proofs since such proofs deal with ranges and domain of unbounded operators to obtain the generator of a $C_0-$semigroup. On other side, this principle remains a extremely powerful tool often encountered in situations  where various bounds are being developed. 

\begin{remark}
For a discussion of consistency in the terms of the resolvent see~\cite{ItoKa98Trotter}; for a discussion of a possible trade-off between consistency and stability, \cite{ItoKa98Trotter, HuItoKa01aspects}; for a version for inhomogeneous ACPs, \cite{HuItoKa01aspects}\footnote{\cite{Schat99stability} gives an example of a non-stable additive splitting scheme.}.
\end{remark}

\begin{remark}
Compositions of generators of contraction semigroups automatically are stable; alternatively, a method is expected to be stable if the operators involved commute~\cite{Bjor98operator}. This is a common duality in the theory of (multiplicative) splitting methods.
\end{remark}

It should also be obvious that  the aforementioned scheme of the classical proof\footnote{However, it is used very often as a general scheme, especially in an analytic setting.} is not the only possible way to establish stronger versions of the Trotter-Kato formula. This topic will be revisited later. 

\begin{remark}
Let $(R_{s,t})_{0\le s\le t}$ be the solution operator of an ACP with initial value $y_0$ and let $(\hat{R}_{k,n})_{k = \ov{0,n-1}, n\in\mbN}$ be solution operators of some recurrent approximation scheme on the uniform (in time) mesh:
\begin{align*}
\hat{y}_{n,k+1} &= \hat{R}_{k,n} \hat{y}_{n,k}, \quad k = \ov{0,n-1}, n\in\mbN,
\end{align*}
where $\hat{y}_{n,n}, n\in\mbN,$ converge to the original solution $R_{0,t} y_0.$
Setting $ \prod_{j=\ov{n,n-1}}\hat{R}_{j,n}=\Id$ we have
\begin{equation}
\label{eq:4.1}
R_{0,t} y_0 -y_{n,n} = \prod_{k=\ov{0,n-1}} \hat{R}_{k,n} (y_0-\hat{y}_{n,0}) + \sum_{k=\ov{0,n-1}}  \prod_{j=\ov{k+1,n-1}} \hat{R}_{j,n} \left(\hat{R}_{k,n} -R_{\frac{k}{n},\frac{k+1}{n}}\right) R_{0,\frac{k}{n}} y_0.
\end{equation}
That is,  the global error is  the sum of propagated local truncation errors 
\[
\left(\hat{R}_{k,n} -R_{\frac{k}{n},\frac{k+1}{n}}\right) R_{0,\frac{k}{n}} y_0, \quad k=\ov{0,n-1}, n\in\mbN,
\]
and the propagated initial error. \eqref{eq:2} is a partial case of this decomposition. 
\end{remark}

\begin{example}
\label{ex:4.meta}
Revisiting \eqref{eq:2}--\eqref{eq:3} we can formulate the following typical  meta theorem: if $A,B,C$ are generators of contraction semigroups $(T_t)_{t\ge 0},(S_t)_{t\ge 0}, (U_t)_{t\ge 0}$ respectively, and for some $x$ and $t\ge 0$
\begin{equation}
\label{eq:4.3}
\|\left(T_{t/n}S_{t/n}-U_{t/n}\right)x\|\le C n^{-p}\|x\|
\end{equation}
 for some $p, C>0$ then
\[
\left\|\left(\left(T_{t/n}S_{t/n}\right)^n-U_{t}\right)x\right\| \le C n^{-p+1}\|x\|.
\] 
Obviously, the tricky part is to obtain the estimate in~\eqref{eq:4.3}, which is often achieved by ad hoc methods. 
\end{example}

\begin{example}(\cite{Faou09Analysis}) 
\label{ex:faou}
Consider the PDE 
\begin{align*}
\frac{\pt u}{\pt t} &= \frac{1}{2}\Delta u + g(u),  \quad u(x,0) = u_0(x), \\
 & x\in \mbR^m, t\ge 0,
\end{align*}
where $g\in C^2(\mbR)$ has bounded derivatives and satisfies $g(0)=0.$ Let $(T_t)_{t \ge 0}$ be the heat semigroup and 
\[
\frac{dF_t(x)}{dt} =g(F_t(x)), \quad F_0(x) =x.
\]
Then it is possible, by using the It\^o formula and properties of the Wiener process, to show that the local error satisfies
\begin{equation}
\label{eq:4.5}
u(x,t) - T_t \left(F_t(u_0(x)) \right) = - \E \int_0^t d F_s\left( u_{t-s}(x + w(s))\right), 
\end{equation} 
where $w$ is a Wiener process, which can be used to derive the consistency bound
\[
\left\| u(t) - (T_t \circ F_t)u_0 \right\|_{L_\infty(\mbR^m)} \le C t^2 \| \nabla u_0\|_{L_{\infty}(\mbR^m)}.
\]
When combined with energy stability estimates for the original semigroup and Example~\ref{ex:4.meta}, this implies that the Trotter-Kato formula is first order accurate.  Note that  the rate of the convergence is obtained by using probabilistic techniques.
\end{example}

\begin{remark}
\cite{Faou09Analysis} also establishes $L_p-$estimates (in the linear case) and results for the $F_t\circ T_t-$version of the Trotter-Kato formula and for the case when one step approximation is $F_{t/2}\circ T_t  \circ F_{t/2}$ (the Strang splitting).
\end{remark}

For other generalizations of the Trotter-Kato, bounds, further developments, extensions to other types of semigroups and references see, besides the sources listed above,~\cite{Ro93error, CaZag01operator, NeidZag98error_estimate, NeiZag99fractional, IchiTa97error, IchiTa98error, JohnLa00Feynman, Ichi04recent, CaNeidZa01accretive, FaouOsSchratz15analysis, NeidZa20convergence,Za05Trotter, HanOs10dimension,HanOs09dimension, NeidZa99fractional, IchiNeidZa04Trotter, IchiTaTaZa01note,IchiTa01norm, Fa67product,Bjor98operator,Vui10generalization,NeidSteZa18operator, VuiWresZa09general,VuiWresZa08Trotter, NeidZa99Trotter, NeidSteZa18remarks, CaNeidZa02comments,NeidSteZa19Trotter,VuiWres11product,CsoEhrFas21operator}
  (both for the Chernoff product formula and the Trotter-Kato formula).

In particular, the non-autonomous case is studied in \cite{DaFo91measures, IchiTa98error, NeidZa20convergence, HanOs09dimension, BatCsoFarNi11operator, EnNa00One-parameter, Fa67product,Vui10generalization, VuiWresZa09general,VuiWresZa08Trotter,NeidSteZa19Trotter, VuiWres11product}. In particular, \cite{VuiWresZa09general,VuiWresZa08Trotter,VuiWres11product} study time-dependent domains. 

Versions of the (non-autonomous) Trotter-Kato formula in alternative settings are established in~\cite{Te68stabilite, DaFo91measures,EisenHan22variational, ItoKa02evolution, BatCsoFarNi11operator}. In particular, \cite{Te68stabilite,EisenHan22variational} use the framework of Gelfand triples (rigged Hilbert spaces).

For results on nonlinear Trotter-Kato formulae, see references in~\cite[Supplement VIII.8]{ReedSai72methods}.

\begin{remark}
Splitting on the level of semigroups means splitting in time. Numerical schemes also include spatial discretization in practice. For results about combining time splitting and spatial discretization (in particular, a version of the Chernoff formula) see~\cite{ItoKa02evolution,BatCsoNi09operator_splittings,BaCsoFarNi12operator_spatial}. See also Remark~\ref{rem:0.1}.
\end{remark}

\begin{remark}
Revisiting examples for Feller processes, one can see that the action of the corresponding semigroup is defined and well behaved for $f\in C_0(\mbR^n)$ ($L_p(\mbR^n)$) and not necessarily for functions outside of such classes though such function may appear in applications. See~\cite{DoerTeich10semigroup} for one extension of splitting schemes to larger classes of functions with preservation of the rate of the convergence for SPDEs.
\end{remark}

\section{A general formulation of a splitting scheme. Multiplicative and additive splitting methods}

The previous sections contain a rigorous  description of the Trotter-Kato formula for semigroups. Henceforward we care more about ideas and illustrations and no longer aim to give precise formulations as we step into the domain of general operator splitting methods as numerical and theoretical schemes to approximate an original solution without any a~priori expectations about convergence or a unified exposition.  

Abstract initial value problems (IVPs)\footnote{We do not call them abstract Cauchy problem as a discussion of semigroups involved is neither provided nor needed.} with no assumptions on the corresponding operators (e.g. they can be multivalued or discontinuous) are used. 
We mostly follow the exposition in~\cite{GlowOshYin17Splitting, HundsVer07Numerical}. 

Consider the autonomous IVP for a possibly nonlinear $A$
\begin{align*}
\frac{du}{dt} + A (u) = 0, \ u(0) = u_0, \quad t\in[0;T].
\end{align*}
 Assume that 
\[
A =\sum_{k=\ov{1,m}} A_k
\] 
for some  $A_k, k=\ov{1,m}.$

In this section we write $A(u)$ and similar terms in the LHS of an  IVP. Later in the text we will resume writing them in the RHS.

Assume $T$ is fixed. Let $h =\frac{T}{N}$ be a fixed time step for some  integer $N.$ Set $t_n = nh$ so $t_0, \ldots, t_{N}$ form a partition of $[0;T].$ 

The following scheme is a general version of the Trotter-Kato formula and is called the {\it Lie-Trotter splitting scheme}. In particular, it is  a the general version of Theorem~\ref{th:lie.product}:
\[
\lim_{n\to\infty}\Big( \e^{\frac{1}{n}A_1} \cdots \e^{\frac{1}{n}A_m} \Big)^n = \e^{\sum_{j=\ov{1,m}} A_j}.
\] 
In the general setting,  define $v_{j,n}, j=\ov{1,m}, n =\ov{0,N-1}$ as follows: for any $n=\ov{0,N-1}$
\begin{align}
\label{eq:5.0}
& \frac{dv_{j,n}(t)}{dt} + A_j v_{j,n}= 0, \quad t \in [t_n;t_{n+1}],  \quad j=\ov{1,m}, \notag \\
& v_{j,n}(t_n) = v_{j-1,n}(t_{n+1}), \quad j=\ov{2,m},  \notag \\
& v_{1,n}(t_n) = v_{m,n-1}(t_{n}),  \notag \\
& v_{1,0}(0) = u_0.
\end{align}
Then 
\[
 u(t_{n+1}) \approx v_{m,n}(t_{n+1}), \quad n =\ov{0,N-1}.
\]  
In general, one expects the Lie-Trotter splitting to be first order accurate. However, this conclusion is not guaranteed.

Adjustments in the non-autonomous case  are standard (though they may be nontrivial mathematically when treated rigorously): e.g. time is considered as an additional variable $\tau(t)$  with $\frac{d\tau(t)}{dt}=1$ and a new operator $(A,-1)$ is introduced or the explicit approach from Example~\ref{ex:2.non} is used. We will use the second option for all other schemes.

Now we consider 
\[
\frac{du(t)}{dt} + A(t, u(t)), \quad t\in[0;T].
\]
Assume $m=2.$ The next scheme is called the {\it Strang splitting scheme} and is asymmetrical in the way $A_1$ and $A_2$ are treated. For matrices, it is Theorem~\ref{th:strang.product}.
To introduce an abstract formulation, consider $v_{1,n}, v_{2,n}, \wt{v}_{1,n}, n=\ov{0,N-1},$ defined via 
\begin{align*}
& \frac{dv_{1,n}(t)}{dt} + A_1(t, v_{1,n}(t))= 0, \quad t \in \Big[t_n, t_{n} + \frac{h}{2}\Big], \\
& \frac{dv_{2,n}(t)}{dt} + A_2(t, v_{2,n}(t))= 0, \quad t \in [t_n, t_{n+1}], \\
& \frac{\wt{v}_{1,n}(t)}{dt} + A_1(t, \wt{v}_{1,n}(t))= 0, \quad t \in \Big[t_n + \frac{h}{2}, t_{n+1}\Big], \\
& v_{1,n}(t_n) = \wt{v}_{1,n-1}(t_{n}),\  v_{2,n}(t_n) = v_{1,n}\Big(t_{n}+\frac{h}{2}\Big),\  \wt{v}_{1,n}(t_n) = v_{2,n}(t_{n+1}), \\
& \wt{v}_{1,-1}(0) = u_0. 
\end{align*}
Then 
\[
 u(t_{n+1})  \approx \wt{v}_{1,n}(t_{n+1}), \quad n = \ov{0,N-1}.
\]
One does expect the second order in the general case but, as earlier, this is not guaranteed. 

We will not repeat this universal remark on the absence of a universal rate for any given splitting scheme in future.

\begin{remark}
The Strang splitting is quite popular not only due to  it being of the second order, but also because it does not require doubling the number of computations (compared to the Lie splitting), as we have already seen in the finite-dimensional case. Indeed, if the propagator for the semigroup of $A_k$ is $(S^{(k)}_t)_t, k=\ov{1,2},$ we  get\footnote{Here $\prod$ is used instead of $\circ$ to make the expression readable.} 
\begin{align*}
u(T) &\approx \wt{v}_{2,N-1}(t_N) \\
&=  S^{(2)}_{t_{N-1}+h/2,t_{N}} \circ S^{(1)}_{t_{N-1}, t_{N}} \circ \prod_{k =\ov{0,N-2}} \Big[ S^{(2)}_{t_k +h/2, t_{k+1} + h/2} \circ S^{(1)}_{t_{k}, t_{k+1}} \Big] \circ S^{(2)}_{0, h/2} u_0 .
\end{align*} 
\end{remark} 

The Lie-Trotter and Strang splitting schemes are {\it multiplicative splitting schemes}: the formal semigroups for $A_1, \ldots, A_m$ are combined via composition exclusively. 

The following {\it Peaceman-Rachford} and {\it Douglas-Rachford operator splitting schemes} are called {\it additive}. We have $m=2.$ To formulate the underlying idea, consider an ODE 
\[
\frac{du}{dt} = f(t, u(t)), \ u(0) =0.
\]
The solution $u$ can be approximated  by using the explicit forward Euler method:
\[
u_{n+1} = u_n + h f(t_n, u_n), \ u(t_n) \approx u_n,
\]
or by using the implicit backward Euler method 
\[
\wt{u}_{n+1} = \wt{u}_n + h f(t_{n+1}, \wt{u}_{n+1}), \ u(t_n) \approx \wt{u}_n.
\]
For $\frac{du}{dt} + Au(t) = 0$ we formally have 
\begin{align*}
& \frac{u_{n+1} - u_n}{h} + Au_n = 0, \\
& u_{n+1} = (A - h\Id)u_n
\end{align*}
and 
\begin{align*}
& \frac{\wt{u}_{n+1} - \wt{u}_n}{h} + A\wt{u}_{n+1} = 0, \\
& \wt{u}_{n+1} = (A + h\Id)^{-1} \wt{u}_n,
\end{align*}
respectively (cf.~\ref{eq:3.euler.b} and \ref{eq:3.euler.f}). 

The idea of the Peaceman-Rachford splitting is as follows: divide $[t_n,t_{n+1}]$ in half,  run the forward Euler scheme for $A_1$
 and the backward Euler scheme for $A_2$ on $[t_n, t_n+\tfrac{h}{2}]$ and run the same algorithm on other half of the interval, switching the roles of $A_1$ and $A_2.$ Setting $\wt{t}_{n} =  t_n +\frac{h}{2}$ we have
\begin{align*}
& \frac{\wt{u}_{n} - u_n}{\tfrac{h}{2}} + A_1(t_n, u_n) + A_2(\wt{t}_n, \wt{u}_{n}) = 0, \\
& \frac{u_{n+1} - \wt{u}_n}{\tfrac{h}{2}} + A_1(t_{n+1}, u_{n+1}) + A_2(\wt{t}_n, \wt{u}_{n}) = 0, \\
& u(t_{n}) \approx u_{n}.
\end{align*}
The autonomous linear version is 
\begin{equation}
\label{eq:5.pea_rac}
u_{n+1} = \Big(1 + \frac{h}{2}A_1\Big)^{-1} \Big(1 - \frac{h}{2}A_2\Big) \Big(1 + \frac{h}{2}A_2\Big)^{-1} \Big(1 - \frac{h}{2}A_1\Big)u_n.
\end{equation}

The Douglas-Rachford splitting uses a similar idea. For $m=2$ we have 
\begin{align*}
& \frac{\wt{u}_{n} - u_n}{h} + A_1(t_n, u_n) + A_2(t_{n+1}, \wt{u}_{n}) = 0, \\
& \frac{u_{n+1} - u_n}{h} + A_1(t_{n+1}, u_{n+1}) + A_2(t_{n+1}, \wt{u}_{n}) = 0, \\
& u(t_{n}) \approx u_{n}.
\end{align*}
The autonomous linear version is 
\begin{equation}
\label{eq:5.dou_rac}
u_{n+1} = \Big(1 + \frac{h}{2}A_1\Big)^{-1} \Big[ \Id - h A_2 \left( \Id+hA_2 \right)^{-1}  \left(\Id - hA_1  \right)   \Big]u_n.
\end{equation}

\begin{remark}
The Douglas-Rachford splitting can  be extended to the case $m>2$\cite{GlowOshYin17Splitting}.
\end{remark}

The last additive scheme we write down explicitly is a development of the Peaceman-Rachford splitting and is called the {\it fractional $\theta-$scheme}: for $\theta\in (0;\frac{1}{2})$
\begin{align*}
& \frac{\wt{u}_n - u_n}{\theta h} + A_1(\wt{u}_n, t_n + \theta h) + A_2(u_n, t_n) = 0, \\
& \frac{\dt{u}_n - \wt{u}_n}{(1-2\theta) h} + A_1(\wt{u}_n, t_n + \theta h) + A_2(\dt{u}_n, t_n + (1-\theta)h) = 0, \\
& \frac{u_{n+1} - \dt{u}_n}{\theta h} + A_1(u_{n+1}, t_{n+1}) + A_2(\dt{u}_n, t_n + (1-\theta) h) = 0, \\
& u(t_n) \approx u_n.
\end{align*}

These additive splitting scheme  are expected to be first order accurate in general and second order accurate if the operators are ''nice''. 

\begin{remark}
If the operator $A_2$ is computationally problematic (e.g. multivalued), both Peaceman-Rachford and Douglas-Rachford splittings can be modified to exclude $A_2$ from the second subproblem. For instance, solving the first equation for $A_2(\wt{t}_n, \wt{u}_{n})$ we get for  the  Peaceman-Rachford splitting
\begin{align*}
& \frac{\wt{u}_{n} - u_n}{\tfrac{h}{2}} + A_1(t_n, u_n) + A_2(\wt{t}_n, \wt{u}_{n}) = 0, \\
& \frac{u_{n+1} - 2\wt{u}_n + u_n}{\tfrac{h}{2}} + A_1(t_{n+1}, u_{n+1}) - A_1(t_n, u_{n}) = 0.
\end{align*}
\end{remark}

There are other common additive operator splitting schemes such as {\it the Tseng splitting} (e.g.~\cite{SuChoGiMouk21parallel, AlGe19iteration}) and {\it the Davis-Yin splitting} (e.g.~\cite{AraTo22direct, ChenChangLiu20three_operator}). 

In general, the additive operator splitting schemes considered above are  known to be convergent provided the operators involved are monotone (or some of them are).

Let us recall some very basic facts about this last assumption and optimization theory~\cite{BauCom17convex, Zeid90nonlinear_2b} and briefly explain some terminology that is often encountered in the context of operator splitting methods in optimization. Let $B$ be a Banach space. A possibly nonlinear and multivalued operator $A\colon B \mapsto B^*$  is called monotone if 
\[
\Re \langle y_1 - y_2, x_1 -x_2 \rangle \ge 0, \quad y_k \in A(x_k), k =1,2.
\]
Operator monotonicity is one of cornerstones of optimization theory in general and convex optimization in particular. The Browder-Minty theorem states that if $A$ is additionally hemicontinuous and coercive\footnote{Recall that these are also standard assumptions in the theory of SPDEs.}, that is, the mappings 
\[
[0;1] \ni s\mapsto \langle A(x + sy), z \rangle, \quad x,y,z \in B,
\]
are continuous and
\begin{align*}
\inf_{x} \frac{\langle A(x), x\rangle}{\|x\|^2} > 0,
\end{align*}
the equation $Au =f$ has a solution. Moreover, the subdifferential $\pt F$ of a proper closed convex function $F$ is  necessarily a (maximally) monotone operator and 
\[ 
u \in \arg\min F \ \leftrightarrow \ 0 \in \pt F u,
\]
so, roughly speaking,  one can solve $\pt F u =0$ instead of the original optimization problem (that is, one should find any function $u$ in the zero set of $\pt F$). Define a proximal operator 
\[
\mathrm{Prox}_F x = \arg\min_y \left(F(y) + \|x-y\|^2 \right). 
\] 
The operator $\mathrm{Prox}_F$ coincides with the resolvent $(\Id + \pt F)^{-1},$ and 
\[
 0 \in \pt Fu \ \leftrightarrow \ u = \mathrm{Prox}_F u.
\] 
Additionally, $\mathrm{Prox}_F = \frac{1}{2}\Id + \frac{1}{2}C$ where $C$ is 1-Lipschitz. Operators that admit such a decomposition are called firmly non-expansive, and fixed point (proximal) iterations are guaranteed to converge for them. The same relations are present for variational inequalities with convex functions.

One finds numerous applications of operator splitting methods in the context of the theory of monotone and proximal operators (in particular, in the context of fixed point iterations). Such results are often formulated either in the terms of fixed points of proximal operators or in the terms of resolvents and zero sets of monotone operators. Technically, however, zero sets can be sampled by running the same fixed point iteration algorithm, and methods based on Lagrange multipliers are commonly used for both formulations.  

Moreover, popular and effective {\it alternating direction methods of (Lagrange) multipliers (ADDM)} are actually the reformulated Peaceman-Rachford and Douglas-Rachford splittings (or other additive operator splitting schemes).  

\begin{example}(\cite{BauCom17convex}) A so-called forward-backward splitting can be derived as follows: we want to find $x$ such that 
\begin{align*} 
& 0 \in (A_1+A_2)x \ \leftrightarrow \  (\Id - A_1)x \in (\Id + A_2 )x \ \leftrightarrow \ \exists y\colon y = (\Id + A_2)x, y = (\Id - A_1)x,   
\end{align*}
so $x$ is such that 
\[
(\Id + A_2)^{-1} (\Id -A_1)x = x.
\]
  Note that alternatively we can show that $x$ is a fixed point for  $(\Id + \alpha A_2)^{-1} (\Id - \alpha A_1)$ for any $\alpha>0.$ Then one may expect that for sufficiently small $\alpha$ fixed point iterations converge not only for monotone $A_1, A_2$  but also if $A_1$ is only cocoercive:  an operator $A$ is $\alpha-$cocoercive if $\alpha A$ is firmly non-expansive. Cocoercive operators often appear as proper gradients.
\end{example}

Note that the corresponding IVPs for all splitting schemes are usually spatially discretized additionally in practice  so we have high-dimensional systems of possibly nonlinear (differential) equations, and one selects methods to solve these systems separately (from the outer splitting scheme). For instance, starting with a second order parabolic PDE
\[
\frac{\pt u}{\pt t} = \tr (A \Hess u)
\] and replacing only spacial derivatives with the corresponding difference quotients we obtain a spatial high-dimensional ODE on some possibly non-uniform grid (or in a space of basis functions if a finite element scheme is used)
\begin{equation}
\label{eq:5.add.1}
\frac{du_{h}(t)}{dt} = B_h u_{h}, 
\end{equation}
where the matrix $B$ depends on the grid and the size $h$ of the grid explicitly\footnote{If a boundary-value problem is considered, a free term appears in the RHS of the equation.}. If both space and time discretizations are used, we have a high-dimensional (possibly nonlinear) equation. 

Time and space discretizations are treated differently on the larger scale: additive methods usually  have a built-in discretization in time while multiplicative splitting schemes, in contrast, do not usually dictate how the subproblems should be solved. Choosing a concrete method leads to a variation of an initial scheme. For instance, by considering the non-autonomous Lie-Trotter splitting and choosing the 1-step backward Euler method we obtain the {\it Marchuk-Yanenko splitting}\footnote{This version is among the earliest numerical splitting schemes historically.}:
define $v_{j,n}, j=\ov{1,m}, n =\ov{0,N-1}$ 
\begin{align*}
& \frac{v_{j,n}-v_{j-1,n}}{h} + A_j(t_{n+1}, v_{j,n})= 0, \quad j = \ov{1,m}, \notag \\
& v_{0,n} = v_{m,n-1}, \notag \\
& v_{m,-1} = u_0, \notag \\
& u(t_{n+1}) \approx v_{m,n}, \quad n=\ov{1,N}.
\end{align*}
In the autonomous case we have an alternative version of the Trotter-Kato theorem
\begin{equation}
\label{eq:5.mar_yan}
v_{m,n+1} = \left( 1+hA_m\right)^{-1} \cdots \left(1+hA_1 \right)^{-1} v_{m,n}.
\end{equation}
Here the outer and inner optimization share the same time step. Alternatively, one can use a Runge-Kutta type method with the step that is smaller than the time step of the splitting scheme etc.

Multiplicative methods that feature a fixed one-step internal optimization step (such as the Marchuk-Yanenko splitting) can be called alternatively {\it locally one-dimensional methods (LOD)} in a generalized sense~\cite{HundsVer07Numerical}, while the Peaceman-Rachford and Douglas-Rachford splittings are known as {\it alternating direction implicit methods (ADI)}.

Multiplicative splitting me\-thods are a subclass of {\it exponential operator splitting methods} (also called simply {\it higher order splitting methods}) -- methods that can be formally represented as linear combinations of formal semigroups\footnote{By saying ''formal'', we only mean that the actual theory of semigroups may not be even invoked in a particular scenario, replaced e.g. with the theory of Lie groups.} of the form\footnote{Here $\prod$ is used instead of $\circ$ to make the expression readable.}
\begin{equation}
\label{eq:5.5}
\sum_{k} a_k \prod\limits_{i} \Big[ \prod_{j=\ov{1,m}} \e^{b_{k,i,j} h A_j} \Big],
\end{equation}
the numbers $\{a_k\}, \{b_{k,i, j}\}$ not necessarily being non-negative or even real. Such methods are often used to construct higher order schemes. 
However, the term ''exponential splitting'' sometimes simply means  the standard Lie-Trotter splitting.

\begin{example}
The autonomous Strang splitting\footnote{We use $(\e^tA)_{t\ge 0}$ to denote the semigroup with a generator $A$ in the context of exponential splitting methods.} 
\begin{equation}
\label{eq:5.6}
R_h=\e^{\frac{h}{2}A_2} \e^{hA_1} \e^{\frac{h}{2}A_2}
\end{equation}
is an exponential splitting scheme. 
\end{example}

\begin{example}
For matrices, the scheme 
\[
\frac{1}{2}\left(\e^{hA_1}\e^{hA_2} + \e^{hA_2}\e^{hA_1} \right) 
\]
is second order accurate. 
\end{example}

\begin{example}\label{ex:20}(\cite{HundsVer07Numerical}) Using~\eqref{eq:5.6}, one can define 
\[
R_{\theta h} R_{(1-2\theta)h} R_{\theta h},
\]
where $\theta >0, 1 - 2\theta <0,$ which is formally of the forth order and requires solving a PDE backward in time. 
\end{example}

Special constructive operator splitting schemes that build on simpler ones (usually multiplicative)  are often called {\it hybrid, weighted, additive} (in the alternative meaning) and {\it iterative} and are briefly treated  in the context of exponential splitting methods later in the text.


{\it Operator splitting methods encountered in the probabilistic setting are mostly multiplicative.} 

The above is a modern classification of operator splitting methods presented in~\cite{GlowOshYin17Splitting}. However, multiplicative and additive methods historically have different origins and different societies of contributors; a specific operator splitting scheme or a particular version of a general scheme (e.g. the St\"ormet-Verlet integration method) may have been developed or rediscovered as a narrow ad hoc method for a particular problem, so the nomenclature is rather diverse and does not necessarily possess integrity a researcher interested in historical accuracy would desire\footnote{even though not to the degree of  W.~Rudin's ironic remark~\cite{Ru91functional}: ''Thus it appears that \v{C}ech proved the Tychonoff theorem, whereas Tychonoff found the \v{C}ech compactification -- a good illustration of the historical reliability of mathematical nomenclature.''}.



The following  short interlude collects some alternative names. 

We start with the Lie-Trotter splitting. A standard and historically justified alternative name is {\it the fractional step method} (or {\it the method of fractional steps}), which can also include other similar schemes.  Other possibilities revolve around (sur)names of S.~Lie, T.E.~Trotter and T.~Kato in different combinations, one of most common being {\it the Lie splitting}~\cite{GlowOshYin17Splitting}. Analysts often speak about {\it the Trotter-Kato formula} only. Probabilists may use {\it ''splitting-up''}. Physicists may speak about {\it the Suzuki-Trotter expansion/decomposition}. {\it Sequential splitting} is another alternative. 

The Strang splitting can be called {\it the Marchuk(-Strang) splitting}. 

The autonomous Peaceman-Rachford and Marchuk-Yanenko splittings  can be written  in the terms of revolvents. Thus {\it the sequential splitting} of~\eqref{eq:5.0}  and {\it the resolvent splitting} of e.g.~\eqref{eq:5.pea_rac}, \eqref{eq:5.dou_rac} or~\eqref{eq:5.mar_yan} can be considered separately. Another example is the Ryu splitting (e.g.~\cite{ArCamTam21strengthened, MaTam23resolvent}).

{\it Additive, sequential, weighted and iterative splitting schemes} may mean different things for different authors and we refer to Section~\ref{sec:exponential} for examples. 
 
Particular versions and algorithms include {\it the St\"ormet-Verlet and leapfrog methods (see Section~10), time-evolving block decimation~\cite{UrSol16parallel, HaHaMcCu20hybrid}, split Hamiltonian\footnote{in particular, split Hamiltonian MC}  methods~\cite{CaSanzShaw22split, ShanLanJohnNeal14split},  split-step Fourier~\cite{AleZaShei05split-step, TaSmithWolf02analysis}, mapping method~\cite{Wis82origin,Mal94mapping}, gradient-projection~\cite{CuiShanb16analysis, FiNoWright07gradient}, split Bregman methods~\cite{GolBreOsher10geometric, OberOsherTaTsai11numerical}, primal-dual splitting~\cite{BotCsetHend14recent, ConKiContHi23proximal}, Rosenbrock's approximate matrix factorization}~\cite{Botver03new, BeckGonPeWei14comparison}, the aforementioned {\it Tseng, Davis-Lin, Ryu splittings}  and others.

\begin{remark}
Intermediate values are approximations of the real solution for additive splitting methods. It is not true for multiplicative schemes.
\end{remark}

\begin{remark}
Additive splitting methods can be seen as a particular instance of implicit explicit mixed methods (IMEX).
\end{remark}

\begin{remark}
Note that the matrix $B_h$ in~\eqref{eq:5.add.1} contains $h^{-2}$ in the case of second derivatives and thus introduces stiffness into the system even if it was not present originally. In particular, implicit methods are advised for solving such ODEs. 
\end{remark}

\section{A naive probabilistic example. An error of a splitting scheme and the sources of it}

This section consider a basic idea behind an actual implementation of a splitting scheme for a parabolic PDE, and the usage of MC for it. 

Since  an additional layer of  spacial discretization is almost always present, the total error of an operator splitting scheme is composed of a theoretical error of the splitting procedure itself, the error of a space(-time) discretization procedure used to obtain  finite-dimensional formulations of subproblems, and  the ordinary calculational error of a chosen method and the corresponding  program implementation used to solve these finite-dimensional problems.  These errors are not additive, obviously, though one can try to obtain an additive bound. On occasion, stability of the final scheme may vary greatly between different choices of internal sub-routines. 

 Discretization here may involve choosing a mesh,  basis expansions, interpolation between points of a space(-time) grid, choosing a specific formula for the first/second/higher derivative at a boundary, variable grid steps, finite differences and finite volumes in general, using volume and mass preservation formulations etc.  Alternatively, if a probabilistic interpretation is available, one may apply MC methods, and  the upper bound of the error becomes, roughly speaking, 
\begin{equation*}
\label{eq:5.4}
E_{splitting} + E_{MC} + E_{sim},
\end{equation*}
where $E_{MC}$ is a statistical discrepancy (due to an insufficient number of simulations), estimated e.g. with using the Berry-Esseen, Bikelis or concentration inequalities~\cite{GraTa13stochastic}, and $E_{sim}$ is the error of the actual simulation procedure for the corresponding random process. 


Let us consider a basic and naive example to explain the situation. Consider an ordinary PDE  with a inhomogeneous term
\begin{align}
\label{eq:8.4}
\frac{\pt u}{\pt t} &= \frac{1}{2} \sum_{k,j =\ov{1,m}} a_{k,j}(x) \frac{\pt^2 u}{\pt x_k \pt x_j} +  \sum_{k=\ov{1,m}} b_k(x) \frac{\pt u}{\pt x_k} + f(x,u)u, \notag \\
u(x,0) &= u_0(x), \quad x\in \mbR^m, 
\end{align}
 and the Lie-Trotter splitting scheme that separates the ODE 
\begin{equation}
\label{eq:8.5}
\frac{dS_t(x)}{dt} = f(x,S_t(x)), \quad S_0(x) = x,
\end{equation}
and the semigroup $(T_t)_{t \ge 0}$ of the Markov process $\xi$ with
\begin{align*}
d\xi_{xk}(t) &= \sum_{j=\ov{1,m}} \sigma_{k,j}(\xi_x(t)) dw_j(t) + b_k(\xi_x(t))dt, \quad k=\ov{1,m}, \notag \\
\xi_x(0) &= x, 
\end{align*}
where we assume that $\|a_{j,k}\|_{j,k=\ov{1,m}} = \sigma \sigma^*$, the square root $\sigma$ is Lipschitz continuous\footnote{A sufficient condition is the non-degeneracy and Lipschitz continuity of $\|a_{j,k}\|_{j,k=\ov{1,m}}.$}, $b_k, k=\ov{1,m},$ are bounded,  $w_j, j=\ov{1,m},$ are independent Wiener processes, and the function $f \in C(\mbR^{2m})$ satisfies suitable growth conditions. Both the process $\xi$ and the semigroup $(S_t)_{t\ge 0}$ are well defined. 

Suppose that the time step of the splitting scheme is $h.$ Then the exact expression for one step of  two possible Lie-Trotter splittings are
\begin{align}
\label{eq:8.7}
\left(T_h \circ S_h\right) u_0(x) &= \E S_h(u_0(\xi_x(h))), \notag \\
\left( S_h \circ T_h \right) u_0(x) &= S_h \left( \E u_0(\xi_x(h))\right). 
\end{align}

Let $\{{S}^{\ve}_t(x)\mid x\in \mbR^d, t \ge 0\}, \ve >0,$ be a family of numerical approximations of the flow $(S_t)_{t\ge 0}$ whose nature may be arbitrary and such that 
\[
{S}^{\ve}_t \to S_t, \quad \ve\to 0,
\]
in some sense. Let $N$ be the number of MC simulations at each step of the splitting scheme and let $\xi_{x,n,1}, \ldots, \xi_{x,n,N}$ be such simulations of $\xi_x(h)$  at time step $n.$ Then one step of the splitting procedure that uses MC to solve~\eqref{eq:8.4} is 
\begin{align}
\label{eq:8.9}
\frac{1}{N} \sum_{k=\ov{1,N}} {S}^\ve_h (u_0(\xi_{x,n,k})) \notag \\
 {S}^\ve_h \Big( \frac{1}{N} \sum_{k=\ov{1,N}}  u_0(\xi_{\cdot,n,k}) \Big)(x),
\end{align}
respectively. The error in these numerical schemes comes from a number of  sources: 
\begin{enumerate}
\item the semigroup $\wt{S}^\ve$ does not recover the original flow unless the ODE~\eqref{eq:8.5} is solved explicitly,
\item the number $N$ of simulations is finite, 
\item the simulations themselves are not perfect:  a random number generator produces pseudorandom numbers actually,  the value of $\xi_{x,n,k}$ is simulated by using smaller that $h$ time steps and there is a limit how small they can get, higher order simulation schemes  require either calculating iterated stochastic integrals reliably or providing families of random variables with the same moments up to a given order etc.,
\item the splitting itself.
\end{enumerate}

However, MC can be unavailable for a specific PDE and even if it is applicable, this method is rarely the first choice for PDEs in practice. 

Note that one cannot really simulate a discretized $u$ on the whole $\mbR^m$ straightforwardly. For instance, standard finite difference schemes should be properly modified  for unbounded domains (e.g. by introducing artificial boundaries or by transforming the original PDE into one on a bounded domain). 

On the other hand, ''real-world'' engineering PDEs  are usually stated for bounded domains automatically.

\begin{remark}
Though MC is significantly less sensitive to the problem of unbounded domains,  the presence of a boundary can cause problems. One can think about using the classical Euler-Maruyama scheme to  simulate a one-dimensional SDE whose solution always stays positive: it is impossible to keep the approximations positive since they use symmetrical random increments, so one has to consider corrected or suitable implicit schemes etc. to remedy this\footnote{See e.g.~\cite{Alfon05discretization, DeNeuSzpruch12Euler} for results and references.}. 
\end{remark}

Consider now \eqref{eq:8.4} in a bounded domain $D\subset \mbR^m$ with BC
\[
u(x,t) = u_b(x), \quad x \in \pt D.
\]
The boundary $\pt D$ is assumed to be piecewise Lipschitz, as is customary. 

We can apply the Feynman-Kac formula to get the analogs of~\eqref{eq:8.7}. 
 $f$ is now defined only on $D\times \mbR,$ so one could try to consider a decomposition (cf. Section 9)
\begin{equation}
\label{eq:8.8}
u_b = \wt{u}_b + \wh{u}_b,
\end{equation}
and a family of degenerate first order PDEs\footnote{For instance, one can add small noise to define this PDE properly.}  to define $\wt{S}_t\circ g =: \wt{S}_t^g$ 
\begin{align*}
&\frac{\pt \wt{S}_t^g(x)}{\pt t} = f(x,\wt{S}^g_t(x)), \\
&\wt{S}_0^g(x) = g(x), \quad  x \in D, \\
&\wt{S}_t^g(x) = \wt{u}_b(x), \quad x \in \pt D,
\end{align*}
instead of~\eqref{eq:8.5}. Since the BC is a function of $x$ only, an alternative is setting $\wh{u}_b=\wt{u}_b=u_b.$    

\begin{remark}
However, references in Section 9 show that such naive approaches to BC are not proper ones. 
\end{remark}

Then \eqref{eq:8.7} transforms into
\begin{align}
\label{eq:8.10}
\left(T_h \circ \wt{S}_h\right) u_0(x) &= \E \wt{S}_h \circ u_0(\xi_x(h)) \1 \left[\tau > h\right] + \E \wh{u}_b\left(\xi_x(\tau)\right) \1\left[\tau \le h\right], \notag \\
\left( \wt{S}_h \circ T_h \right) u_0(x) &= \wt{S}_h \circ \Big( \E u_0(\xi_\cdot(h)) \1 \left[\tau > h\right] + \E \wh{u}_b\left(\xi_\cdot(\tau)\right) \1\left[\tau \le h\right] \Big)(x),
\end{align}
where $\tau$ is the time when $\xi_x$ hits $\pt D.$ 

\begin{remark}
Even though we consider only the Dirichlet BC, it is clear  that the decomposition~\eqref{eq:8.8} cannot be arbitrary even if $f$ can accommodate it. E.g. if $\sigma=0$ on some $C\subset \pt D,$ this set $C$ may be unreachable for $\xi$ (so BC does not make sense here at all) or the process $\xi$ cannot be reflected there (so the Neumann BC cannot be imposed) etc. Recalling that no classification of boundary points of a diffusion exists in dimensions higher than or equal to $2,$ we can see this as a particular difficulty of the probabilistic origin.
\end{remark}

 
Assume also that we are given a family of possibly non-uniform grids (tessellations) $G_\ve, \ve >0,$ on $D$ whose sizes converge to $0$ as $\ve\to 0.$  We have the following ODE for the finite-dimensional spatially discretized version of $(T_t)_{t \ge 0}$
\begin{align}
\label{eq:8.2}
\frac{d{v}_\ve(t, a)}{dt} &=  B_\ve({v}_\ve(t,a)) {v}_\ve(t,a) + g_\ve({v}_\ve(t,a)), \notag \\
 v_\ve(0,a) &= v_{\ve 0}(a), \notag  \\
 &t \in [0;h], \ a\in G_\ve,
\end{align}
where the term $g_\ve$ appears due to boundary corrections when calculating derivatives and explicitly depends on $\wh{u}_b,$ and $v_{\ve, 0}$ is a discretized version of the IC (which changes at each step of the splitting scheme). Let ${v}_{\ve, n}$ denote finite-dimensional approximations of  $(T_t)_{t \ge 0}$ on the grid $G_\ve$ at time step $n.$   One can use the standard  $\theta-$scheme\footnote{This is the Crank-Nicolson scheme if $\theta=\frac{1}{2}.$} with time step $h$ for solving ODEs to get
\[
{v}_{\ve,n+1} = {v}_{\ve,n} + \theta h  \left[ B_\ve({v}_{\ve,n+1}) {v}_{\ve,n+1} + g_\ve({v}_{\ve,n+1}) \right] + (1-\theta) h \left[ B_\ve({v}_{\ve,n}) {v}_{\ve,n} + g_\ve({v}_{\ve,n})\right], 
\]
which should be solved for ${v}_{\ve,n+1}$. As it has been remarked already, the matrix $B_\ve$ depends on $G_\ve$ and contains the second power of the grid size of $G_\ve.$ Thus it needs implicit and sufficiently stable methods whose step sizes get smaller as $\ve\to 0$ and thus  limit  the choice of the grid in practice. 
 This introduces an error the value of which depends on $G_\ve,$ the domain $D,$  properties of $\pt D$  etc. Such schemes are run at every step of the splitting procedure. 

However, let us choose for simplicity the forward Euler method with time step $h$ so we have instead just a linear matrix-valued equation
\begin{align}
\label{eq:8.3}
& \frac{{v}_{\ve,n+1} - {v}_{\ve,n}}{h} =  B_\ve({v}_{\ve,n}) {v}_{\ve,n} + g_\ve({v}_{\ve,n}), \notag \\
& {v}_{\ve,n+1} = T^\ve_h ({v}_{\ve,n}), \notag \\
& T^\ve_h(x) = \Big(\Id + h B_\ve(x) \Big)x + h g_\ve(x), \quad x \mathrm{\ on\ } G_\ve.
\end{align}

We also need $\{\wt{S}^{\ve}_t(x)\mid x\in \mbR^d, t \ge 0\}, \ve >0,$ which is a family of  numerical approximations of the ''flow'' $(\wt{S}_t)_{t\ge 0}$ whose nature may again be arbitrary and such that 
\[
\wt{S}^{\ve}_t \to \wt{S}_t, \quad \ve\to 0,
\]
in an appropriate  sense. These approximations are also finite-dimensional and live on its own grids so $\wt{S}^\ve_h(x)$ may not be defined for an arbitrary $x\in D.$ But $T^\ve_h(x)$ also makes sense only if $x\in G_\ve.$
Therefore to combine the operators $\wt{S}^\ve_h$ and $T^\ve_h$ from~\eqref{eq:8.3} we need to introduce interpolation procedures $I^\ve_h, J^\ve_h, \ve >0,$ such that $\wt{S}^\ve_h \circ I^\ve_h \circ T^\ve_h,  T^\ve_h \circ J^\ve_h \circ \wt{S}^\ve_h$
are well defined. The simplest example is a linear interpolation on a uniform grid. 

Then we get instead of~\eqref{eq:8.10}
\begin{align}
\label{eq:8.12}
& T^\ve_h \circ J^\ve_h \circ \wt{S}^\ve_h, \notag \\
& \wt{S}^\ve_h \circ I^\ve_h \circ T^\ve_h.
\end{align}

It is left to write down the corresponding versions of MC
\begin{align}
\label{eq:8.11}
& \frac{1}{N} \sum_{k=\ov{1,N}} J^\ve \circ \wt{S}^\ve_h \circ I^\ve \circ u^\ve_0(\xi_{\cdot,n,k})(x), \notag \\
& J^\ve \circ\wt {S}^\ve_h \circ I^\ve \Big( \frac{1}{N} \sum_{k=\ov{1,N}}  u^\ve_0(\xi_{\cdot,n,k}) \Big)(x),
\end{align}
where $u^\ve_0$ is now defined on $G_\ve$ only.

Thus we get additional sources of error:
\begin{enumerate}
\item the decomposition of BC and the corresponding corrections (if any);
\item solving the ODE part (and the corresponding PDEs);
\item interpolation between points of a grid;
\item space discretization, properties of the grids, the choice of difference quotients etc.;
\item internal time discretization (and associated internal errors).
\end{enumerate}

We will continue the discussion of how BC affects splitting later.

In practice, the grid (that is, $\ve$) is fixed and choosing a proper $\ve$ means choosing a sufficiently numerically efficient  spatial discretization. However, different subproblems (a diffusion and an ODE in our case) may use different spatial discretizations on different grids. For instance, assume that a dimension splitting is applied and the second order operator is decomposed itself so each subproblem is a parabolic IVP in a subspace.      

All these schemes --  \eqref{eq:8.7}, \eqref{eq:8.9}, \eqref{eq:8.10}, \eqref{eq:8.11}, \eqref{eq:8.12} -- are different interpretations of the same splitting procedure. Regardless of the path, the end result shows the discrepancy between a theoretical analysis of our abstract splitting scheme and a practical implementation of it, and the same discrepancy is present between any abstract splitting method and its implementation. For instance, if a Faedo-Galerkin scheme is used, (2) is replaced with properties of the scheme, bases and spaces involved etc. 

\begin{remark}
If an additive splitting scheme is used to solve an optimization problem, one still calculates only a discretized approximation of the optimal function etc.
\end{remark}

Actually, one can initially use approximate discretized semigroups to formulate a splitting scheme. This is how matrix-valued splitting schemes appear. \cite{Pazy12Semigroups} can be consulted for such a formulation in the case of standard approximation schemes and corresponding applications of the Trotter-Kato theorems and e.g.~\cite{ItoKa02evolution,BatCsoNi09operator_splittings, BaCsoFarNi12operator_spatial, EisenHan22variational, KoLu17numerical} can be consulted for abstract results for splitting methods in this setting.

Again, a particular splitting scheme can be selected first in practice, and  optimization subroutines may follow naturally. Then different parts of the pipeline can be tuned separately. 

\begin{remark}
One important observation is that spatially discretized operators are often just bounded linear operators, and finite-dimensional results about the convergence are needed for them.  
\end{remark}

\begin{remark}
In the case of a general PDE and a general method (e.g. the finite element method), an additional matrix may appear in~\eqref{eq:8.2} before the time derivative~\cite{HundsVer07Numerical}, so even the ''explicit'' forward Euler method  of~\eqref{eq:8.3} would use matrix inversion, which  requires some additional effort.
\end{remark}


\begin{remark}
A separation of linear  and nonlinear parts for a Cauchy problem (e.g. for advection-diffusion-reaction equations) is often a viable strategy (see e.g. \cite{HundsVer07Numerical, Schatz02numerical, GlowOshYin17Splitting, DesDuDuLouMa11adaptive, DesDuDuLouMa14analysis} and references therein) in general due to technical and numerical considerations and may even lead to theoretical results by directly exploiting properties of the flow of the ODE that corresponds to the nonlinear part (e.g.~\cite{Faou09Analysis} and Example~\ref{ex:faou}). For instance, this ODE may be solvable explicitly.
\end{remark}

\begin{remark}
An important and frequent task is to find a steady-state/equilibrium solution with $\frac{\pt u}{\pt t} = 0.$ Such solutions often represent the long-term behavior of a physical system. In general, operator splitting methods act on a fixed time interval, and it is advised to use other methods to seek such solutions or to provide a modified version of a splitting scheme~\cite{GlowOshYin17Splitting, Bo23analysis, HundsVer07Numerical}. 
\end{remark}

\section{Some remarks on the relation between ACPs and practice}
We used ACPs to introduce the Trotter-Kato product formula (and thus the Lie-Trotter splitting scheme) rigorously. They provide a well-established framework which can be extremely useful in theoretical studies or applications. However, even multiplicative operator splitting methods are not limited to this framework due to a number of reasons, some of them being limitations of the semigroup approach. 
  Even though a semigroup representation is available for a huge number of PDEs including those involving nonlinear and multivalued operators, it may be  preferable to use the weak/variational formulation to describe optimization problems, PDEs and especially partial differential inequalities. 

Indeed, the possibility to use finite difference and finite element methods (Faedo-Galerkin schemes), gradient descent-based procedures, reformulations in the terms of dual problems of convex programming, parallelization and so on often -- but not necessarily -- pairs well with the variational approach. 



These two intertwined frameworks compliment each other. For instance, a variational reformulation of the Trotter-Kato theorem can be found in~\cite{ItoKa02evolution}. One can recall two different formulations of SPDEs as an example of the interplay between these two competing frameworks~\cite{LiuRock15Stochastic, GaMand11Stochastic}. 


\begin{remark}
It is sometimes beneficial to introduce another types  of solutions (e.g. viscosity solutions, entropy solutions), the  Navier-Stokes equation being an example where even the discussion of well-posedness is nontrivial.
\end{remark}


It also must be noted that it may be of crucial interest for a researcher to have a decent and not necessarily optimal simulation procedure for a specific applied problem.
Moreover, it may happen also that only a limited set of (cost) effective, sufficiently fast or (programmatically) implemented solvers is available, which dictates a concrete (and possibly theoretically subpar and non-optimal) decomposition of an original problem.

In the end, even though the rate of convergence is often of utmost importance in applications, it may nevertheless happen that the precise rate of  convergence is unknown or that even the conclusion about the convergence of an operator splitting scheme is missing in a particular case entirely -- and yet the sheer possibility to develop a sufficiently reliable and practically implementable working simulation procedure can suffice to justify the usage of this scheme. For instance, it may be sufficient to have only the convergence result for the scheme itself.

To conclude, we have a difference between theoretical difficulties and engineering challenges\footnote{Such an extreme pathology as the Lebesgue thorn of the potential theory cannot  appear in real life; theoretical technical advances on irregular and regular boundary points (in the sense of Wiener) for domains with corners are rarely of interest for pure engineers etc.} and a gap between available rigorous theoretical results and practical needs. However, applications  is exactly the role  operator splitting methods excel in besides theoretical studies. In fact, the development of additive operator splitting methods and semi-discretized multiplicative ones  was motivated and driven by practical studies. 



\section{The rate of convergence. Some remarks about error representations in  the deterministic case}
It is nontrivial to find the accuracy of a splitting scheme. We have already discussed the approach that combines the local error and the stability estimates. The reader can consult Sections~4-5 for many references which typically use the theory of semigroups. To obtain the rate of convergence, nontrivial additional assumptions  should be added. 

Some results on the speed of convergence (and occasionally on directly related topics such as stability preservation) that have not been mentioned yet earlier and are at least to some degree outside the analytic (semigroup-based) setting of the classical Kato-Trotter formula are \cite{EiSchna18error,EiJahnSchnau19error,FaouOsSchratz15analysis, FaGei07iterative, IchiTa97estimate, IchiTa00estimate_2,SheAza20splitting, EinOs15overcoming_1, EinOs16overcoming_2, OsSchratz13stability, HanOs16high, OsSchratz13stability, KuhWa01commutator, HanOs08dimension, Ta97error, AzuIchi08note, IchiTa97estimate, IchiTa06exponential,IchiTa98error, DiaSchat96commutateurs, DiaSchatz9estimations, Sheng94Global, GaLaMon19unbalanced, Gei08fourth, HolKarlLie00operator,OsSchratz12error, KirWads02solution,Gei15recent,La16convergence,HaSu05finding, Su91general, Su92general,Su95hybrid,Su93improved,Su91general,FaHa07consistency,BlaCaMu08splitting, BlaCaCharMu13optimized, BatCsoEn12Stability,BatCsoNi09operator_splittings,BaCsoFarNi12operator_spatial, EinMoOs18efficient, KochLu08variational, HanHe17additive, CsoEhrFas21operator}. See also references later in this section. References given in Section 1 (\cite{GlowOshYin17Splitting, GloLeTa89Augmented, Glo03Finite, Ya71Method, RyuYin23large_scale, BauCom17convex, Mar82methods, Mar88splitting, Mar86splitting} and others) can be used as a starting point for anyone interested in splitting methods in numerical methods (especially in the field of fluid mechanics and other physical applications and/or variational methods), and the corresponding literature is extremely rich. 

\begin{remark}
In particular: \cite{IchiTa97estimate, IchiTa00estimate_2} study the difference between the Feynman-Kac formula and the Shr\"{o}dinger operator using probabilistic techniques including the theory of Lev\'y processes (see also references therein); \cite{Ta97error, AzuIchi08note, IchiTa97estimate,IchiTa06exponential} obtain estimates on the difference of the corresponding kernels, also via probabilistic techniques. 
\end{remark}

As illustrations, we include the following results about applications (and occasionally about implementation) of splitting methods:\footnote{We do not include pure simulations. We try to give a wide range of applications or methods and sometimes prefer recent publications since they can easy lead to older results in case a reader is interested. A specific important equation may be understood as an ''application''.}~\cite{DuMaDes11parareal, Schatz02numerical, JahnAl10efficient, BouCar14splitting, ShiLiZhaoZou14new, PengZhaoFeng17operator, Gei08fourth, Gei11computing, Gei11iterative_book, KarlLieNatNordDah01operator, HolKarlLie00operator_2, GiRe19operator, Ka05accurate,CaiFangChenSun23fast, AltZi23secondArx,AltKoZi23bulk, BaCsoFar13operator,GonHerPe14Rosenbrok, HeLawDraPet14local, MarTeichWol16parabolic, HanKraOs12second,LiuMasRi23convergence, ChoiChoiKoh23convergence, LiuWangWang22second_order, DuanChenLiuWang22second_order, LiuWangWangWise22convergence, AuKaKochThal17convergence, Ze22analysis, KochLu11variational, KoLu17numerical}. 

In general, the rate of the convergence can be arbitrary (or a bound may not be available at all). E.g.~\cite{LiuWangWang22second_order} establishes only the well-posedness of a scheme and its energy-preservation properties (cf.~\cite{NeidZa99Trotter}). 

\begin{remark}\label{rem:0.1}
\cite{LiuMasRi23convergence, DuanChenLiuWang22second_order, LiuWangWangWise22convergence, AuKaKochThal17convergence} provide a bound of the speed of the convergence that combines sizes of time and space discretizations.
\end{remark}

The following 3 examples illustrate the situation with the Trotter-Kato formula.


\begin{example}\label{ex:0.2}(\cite{Za05Trotter})
Let $A_1, A_2, A=A_1+A_2$ be non-negative self-adjoint operators. If 
\[
D((A_1+A_2)^\alpha) \subset D(A_1^\alpha) \cap D(A_2^\alpha)
\]
for some $\alpha \in (\tfrac{1}{2};1),$ then
\[
\left\|\left(\e^{-\frac{1}{n}A_1}\e^{-\frac{1}{n}A_2}\right)^n - \e^{-A} \right\| \le\frac{C}{n^{2\alpha-1}}.
\]
This covers the example of two Laplacians on a bounded domain with a smooth boundary, one with a Dirichlet BC and the other with a Neumann BC ($\alpha \in (\frac{1}{2};\tfrac{3}{4})$).
\end{example}

\begin{example}\label{ex:0.1}(\cite{IchiTaTaZa01note,IchiTa01norm}, cf.~\cite{Zag22notes})
Let $A_1, A_2$ be non-negative self-adjoint operators such that $A=A_1+A_2$ is self-adjoint on $D(A_1)\cap D(A_2).$ Then
\[
\left\| \left(\e^{-\frac{1}{n}A_1}\e^{-\frac{1}{n}A_2}\right)^n - \e^{-A} \right\| \le\frac{C}{n}.
\]
This assumption on $A$ is stronger than in the previous example.
\end{example}

\begin{example}(\cite{HanOs08dimension}, see also~\cite{OsSchratz13stability, FaouOsSchratz15analysis, EinOs15overcoming_1, EinOs16overcoming_2})
Consider an ACP on a Hilbert space
\[
\frac{du}{dt} = (A_1+A_2)u,
\]
where $A_1,A_2$ and $A=A_1+A_2$ generate $C_0-$semigroups. If
\begin{align*}
& D(A^2) \subset D(A_1 A_2), \\
& \| A_1 A_2 (\Id - A)^{-2} \| < \infty,
\end{align*}
then the Lie-Trotter splitting is first order accurate. As an example, the coordinate-wise decomposition of a possibly degenerate parabolic operator with a Dirichlet BC can be considered. 
\end{example}

\begin{remark}
A characteristic feature of Example~\ref{ex:0.2}, in addition to dealing with fractional power of operators, is that both operators are comparable and thus the perturbation theory that is sometimes encountered in the setting of splitting methods cannot be used.
\end{remark}


\begin{remark}
There are some other classes of semigroups (e.g. holomorphic ones) that the rate of the convergence in the Trotter-Kato formula is known for.
\end{remark}

\begin{remark}
See also the survey~\cite{CsoSi21numerical}.
\end{remark}

Though  most references above and from the previous sections are heavily analytic, it is not obligatory to follow this route and a variety of methods only vaguely related to the semigroup theory\footnote{which is itself a wide field} or  totally different approaches are used.

The next example, besides its direct purpose, also introduces basics of the variational formulation of IVPs.

\begin{example}\label{ex:0.3}(\cite{EisenHan22variational}, cf.~\cite{Te68stabilite}; also see~\cite{LiuMasRi23convergence, DuanChenLiuWang22second_order, LiuWangWangWise22convergence, AuKaKochThal17convergence, HanHe17additive, KoLu17numerical})
A shortened exposition of a general variational framework for  one semi-discrete multiplicative splitting scheme with the backward Euler internal step is as follows.  
Assume that separable reflexive Banach spaces $V_j, j=\ov{0,m},$ with $\cap_{j=\ov{1,m}}V_j = V_0,$ are continuously and densely embedded into a Hilbert space $H$ so one can consider Gelfand triples 
\[
V_j \hookrightarrow H \cong H^* \hookrightarrow V_j^*, \quad j=\ov{0,m}.
\]
The norms $\|\|_{V_0}$ and $\sum_{j=\ov{1,m}} \|\|_{V_j}$ are assumed to be equivalent. For fixed $j,$ let $A_j(t), t\in [0;T],$ be hemicontinuous  coercive  monotone\footnote{see Section~5} operators from $V_j$ to $V_j^*$ such that $\sum_{j=\ov{0,m}} A_j =A_0$ on $[0;T],$ and for some $p>1$
\[
\|A_j(t)v\|_{V^*_j} \le C(1 + \|v\|^{p-1}_{V_j}), \quad v\in V_j, j=\ov{0,m}.
\]
Set $q = \frac{p}{p-1}.$ Let $f_j \in L^{q}((0;T), V^*_j), j=\ov{0,m},$ be such that 
\begin{align*}
& \sum_{j=\ov{1,m}} f_j = f_0, \\ 
& \|f_j\|_{V^*_j} \le \|f_0\|_{V^*_0} \ \ \mathrm{a.e.}
\end{align*}
Then, in particular, each equation
\[
\frac{du_j}{dt} = Au_j + f_j 
\]
has the unique solution $u_j \in L^p((0;T), V_j)$ with $u_j^\prime\in L^q((0;T), V_j^*).$ To introduce time discretization, set for fixed $n\in\mbN$
\begin{align*}
& A_{j,n,k} = A_j\left(\frac{k}{n}\right), \\
& f_{j,n,k} = n\int_{\frac{k}{n}}^{\frac{k+1}{n}} f_j(s)ds, \quad k =\ov{0,n}.
\end{align*}
The splitting schemes is then defined as follows: for any $n\in\mbN$ 
\begin{align*}
& n\left( u_{j,n,k} - v_{n,k-1} \right) = m \left( A_{j,n,k}u_{j,n,k} + f_{j,n,k} \right) \quad \mathrm{in\ } V_j^*, \\
& \quad j=\ov{1,m}, \\
& v_{n,k} = \frac{1}{m} \sum_{j=\ov{1,m}} u_{j,n,k}, \\
& v_{n,k} \approx u_0\Big(\frac{k}{n}\Big), \quad k =\ov{0,n}.
\end{align*}
Piecewise constant and linear prolongations of $\{u_{n,k}\}$ converge to $u_0$ pointwise in $H.$ If the operators $A_j(t)$ are strongly monotone, the convergence also is in  $L^p((0;T), V_0).$
\end{example}

\begin{remark}
If the operators in Example~\ref{ex:0.3} are discretized in space and thus finite-dimensional, the Gelfand triples collapse to Euclidean spaces.
\end{remark}

The general analysis of parabolic problems in domains with arbitrary BC is in development\ and assumptions of general theorems should be checked in every particular case and are not guaranteed to apply. 

Now we present a short illustrative description of one very different method/tool of studying errors of a splitting scheme which has not been represented in the previous references (with a few exceptions).

Recall that  ODEs with constant coefficients is covered by \eqref{eq:2.add.1} and \eqref{eq:2.add.2} where explicit expressions for local errors in the terms of commutators are given. To extend the method towards general ODEs (or a PDE) of the form
\[
\frac{du}{dt} = f_1(t, u) + f_2(t, u),
\]
one can use the ordinary Taylor expansion. The result often includes iterated commutators  and estimating such expressions is a vital part of the approach.

\begin{example}\label{ex:10.1}(\cite{HundsVer07Numerical}, cf. e.g.~\cite{Bjor98operator,FaHa07consistency,SheAza20splitting})
In the autonomous one-dimensional case this idea gives for the Lie-Trotter splitting with time step $h$
\begin{equation}
\label{eq:10.1}
\frac{h^2}{2}[f_1, f_2] + O(h^3),
\end{equation} 
where $[f_2, f_1] = f_1^\prime f_2 - f_2^\prime f_1.$
\end{example}

\begin{example}\label{ex:10.2}(\cite{Yo90construction})
Consider the Strang splitting $S_{2}(\tfrac{h}{2})S_{1}(h)S_{2}(\tfrac{h}{2})$ for
\[
\frac{dS_k(t)}{dt} = f_k(S_k(t)), \quad k =\ov{1,2}.
\]
 The local error is
\begin{equation}
\label{eq:10.2}\frac{h^3}{12}[f_2, [f_1, f_2]] + \frac{h^3}{24} [f_1, [f_2, f_1]] + O(h^4).
\end{equation}

\end{example}

\begin{remark}
The algebraic interpretation of \eqref{eq:10.1} and \eqref{eq:10.2} is discussed in Section~\ref{sec:exponential}. 
\end{remark}

A related and often used simultaneously approach is to derive PDEs for global or local errors. Such PDEs, being treated as IVPs or ODEs in functional spaces, usually lead to expressions for the errors in the terms of the variation-of-constants formula and again involve multiple iterated commutators of differential operators. Alternatively, an error can be given as a difference of  integral representations of the exact solution and the approximation. This formally works  for almost any operators but one still needs to prove the validity of such formulae and obtain bounds for them. 

\begin{example}(\cite{Sheng94Global})
The error of  the Lie-Trotter splitting is
\begin{equation}
\label{eq:add.100.1}
\e^{t(A_1+A_2)} -\e^{tA_1}\e^{tA_2} = \int_0^t \e^{(t-s)(A_1+A_2)} \left[A_2, \e^{sA_1} \right] \e^{sA_2}ds. 
\end{equation}
\end{example}

General examples are \cite{DesSchatz02Strangs, DesMa04operator, DuMaDes11parareal, DesDuLou07local, Des01convergence, HanOs16high, KochNeuThal13error, JahnLu00error, KochLu08variational, AuKochThal12defect_1, AuHofKochThal15defect_3, AuKochThal15defect_local, CharMeThalZhang16improved, DesDuDuLouMa11adaptive, DesThal13Lie, Thal08high}.

\section{Boundary conditions, inhomogeneous ACPs  and the phenomenon of order reduction}
Now we return to the discussion of the decomposition of BC. 

In general,  BC often leads to the phenomenon of order reduction  when the rate of the convergence of a splitting scheme is observed to be lower that expected~\cite{HundsVer07Numerical}: e.g. the order of the Strang splitting is strictly lower than $2$  (it may be only first order accurate in practice~\cite{FaouOsSchratz15analysis}), while the Lie-Trotter splitting may happen to preserve the first order in the same setting.

 Recall that a (semi-)discretized version of a PDE in a domain can be considered as an inhomogeneous ODE where BC is explicitly used to derive the RHS. Thus the following example  shows that a naive and straightforward approach leads to unavoidable additional errors even in the simplest cases both in the context of inhomogeneous PDEs (ODEs) and non-trivial BC (as well as the corresponding numerical schemes for such equations). 

\begin{example}\label{ex:9.1}(\cite{HundsVer07Numerical})
Assume that matrices $A_1,A_2\in M_n$ commute and consider
\begin{align*}
& \frac{du(t)}{dt} = (A_1 + A_2)u(t) + g(t), \\
& u(t) = \e^{t(A_1+A_2)}u_0 + \int_0^t \e^{(t-s)(A_1+A_2)} g(s) ds.
\end{align*}
with $g=g_1 + g_2.$ One step of the Lie-Trotter splitting with step size $h$ is
\begin{align*}
& u_{n} = \e^{h(A_1+A_2)}u_{n-1} + \e^{hA_2} \int_0^h \e^{(h-s)A_1}  g_1(hn+s) ds + \int_0^h \e^{(h-s)A_2}g_2(hn+s)ds.
\end{align*}
Thus the local error is obviously non-zero even for constant $g_1, g_2.$  However, the method becomes exact if we replace $g_1, g_2$ with
\begin{align}
\label{eq:9.2}
\wt{g}_1(hn+s) &=  \e^{-sA_2} g_1(hn+s), \notag \\
\wt{g}_2(hn+s) &= \e^{(h-s)A_1} g_2(hn+s),  \notag \\
& s \in [0;h].
\end{align}
\end{example}

On the other hand, the theoretical analysis of ACPs often incorporates BC into the definitions of the corresponding domains, so they do not appear in the equation directly. Otherwise a corrective additional inhomogeneous term appears. 

\begin{example}(\cite{EinOs16overcoming_2})
Consider the Laplacian in a domain $D$ with IC $u_0$ and BC $f.$ The solution  can be represented as  $v+w$ where $w$ is the harmonic extension of $f$ in $D$ and $v$ is the solution of
\begin{align*}
& \frac{\pt v}{\pt t} = \Delta v + g(v,w),\\
& v(0) = u_0 - w_0 \quad \mathrm{in\ } D, \\
& v = 0 \quad \mathrm{on\ }  \pt D, 
\end{align*}
with the correction $g$\footnote{See \cite{FaouOsSchratz15analysis} for a discussion of the situation when $w$ is problematic to find.}. However, $\frac{\pt v}{\pt t} =g$  do not necessarily satisfy the original BC and may thus be incompatible with the Laplacian (in the sense of domains), which is a problem for the theoretical error analysis.
\end{example}

\begin{remark}
Note that this point of view  can be somewhat different from most typical approaches based on the weak formulation and Galerkin-Faedo approximations where BC is often imposed onto ODEs for basis coefficients or is included via regularization. For a quick summary of the relation between inhomogeneous PDEs, non-trivial BC, the weak formulation of a PDE, Galerkin-Faedo appproximations and other common numerical methods, see e.g. references in~\cite{PeetPeet24treatment}. See also \cite{LuSiWahl92semigroup} for the semigroup approach.
\end{remark}

Example~\ref{ex:9.1} shows that one can try to eliminate this source of errors by correcting boundary terms  at each step of a splitting scheme.  However, \eqref{eq:9.2} contains $\e^{-sA_2},$ which means backward integration in time (solving the initial IVP backward in an abstract setting). This is  not only a problem for implicit methods in general but even such a standard operator as the Laplacian does not permit such integration. Obviously, the same remarks apply to the decomposition of the inhomogeneous part of a PDE in general, even without BC.

\begin{example}\label{ex:9.2}
The solution of a one-dimensional heat equation can be written as a series that  diverges for negative $t.$ This is a standard simple  example of a PDE that cannot be solved backward in time. 
However, the heat semigroup is analytic and thus accepts some complex numbers with positive real part.
\end{example}


\begin{example}(\cite{FaouOsSchratz15analysis})
Dimension splitting for an inhomogeneous two-dimensional parabolic PDE\footnote{in the divergence form, which is typical for such analytical results} in a domain is considered.
 For an inhomogeneous ACP and in the same notation as in Example~\ref{ex:9.1}, the corrected Lie-Trotter splitting  
\[
u_{n+1} = \e^{hA_1}\e^{hA_2}(u_n + hg(t_n)) 
\]
is first order accurate under suitable regularity assumptions on $g.$ The order of accuracy for the corrected Strang splitting
\[
\e^{\frac{h}{2}A_1}\e^{\frac{h}{2}A_2}\Big( \e^{\frac{h}{2}A_2}\e^{\frac{h}{2}A_1}u_n + hg(t_n + \tfrac{h}{2}) \Big) 
\]
is guaranteed to be less than  $2$  unless additional assumptions on the smoothness of $g$ and $u(0)$ are imposed. 
\end{example}

Some other references are \cite{HanOs09dimension,EinOs15overcoming_1, EinOs16overcoming_2,OsSchratz13stability, EinOs18comparison,OsSchratz12error,AltZi23secondArx,AltKoZi23bulk,EinMoOs18efficient,HanKraOs12second}. In particular, \cite{EinOs18comparison} provides a comparison of some methods for the Strang splitting, and \cite{AltZi23secondArx,AltKoZi23bulk}  study BC that directly involve the solutions. In particular, \cite{OsSchratz13stability} develops general results on stability of  splitting schemes in the terms of smoothing properties of analytic semigroups involved (that is, in the terms of fractional powers of the corresponding generators);  domains of various operators that appear in this context should also be suitably  compatible\footnote{Cf. the discussion of the control for domains in the proof of the Trotter-Karo formula in Section 4.}. Typically, the BC and inhomogeneous part of the equation should be sufficiently smooth. A  problem of BC with time derivatives is considered in~\cite{CsoEhrFas21operator, AltZi23secondArx,AltKoZi23bulk, KoLu17numerical}.

\begin{remark}
A very naive -- and obvious -- explanation of some results here is as follows: to obtain higher order stability results, one needs to expand a solution (given e.g. via the variation-of-constants formula or via \eqref{eq:add.100.1}) and estimate the remainder, which yields assumptions on the regularity of the solution. This produces requirements for BC and IC, and it is exactly where fractional powers of closed operators appear. See e.g.~\cite{OsSchratz13stability} for a detailed example. Probabilists also use such techniques in the theory of SPDEs and in the theory of subordinated Feller and sub-Markovian semigroups~\cite{But20Method, Ja02pseudoVol2, Ja02pseudoVol3}\footnote{Subordination is  a random time change by a subordinator.}, for example. In particular, see \cite{But20Method, Ja02pseudoVol2} and other related references in Section~4 for applications of the Trotter-Kato formula in the latter context.
\end{remark}

\begin{remark}
Some detailed numerical splitting schemes can be found e.g. in~\cite{EinMoOs18efficient}.
\end{remark}



\section{General exponential splitting methods. The  Baker-Campbell-Hausdorff  formula}
\label{sec:exponential}

This section discusses the general idea behind exponential splitting methods (recall~\eqref{eq:5.5}) and the  Baker-Campbell-Hausdorff (BCH) formula\footnote{also known as the Baker-Campbell-Hausdorff-Dynkin formula} in particular.

The underlying idea is to expand exponentials in~\eqref{eq:5.5} for small times and to choose coefficients and times of integration 
in such a way that all lower order terms cancel. 

Such task admits a full solution in the special case of ODEs that possess additional geometric properties. 


\begin{example}
Given $A\in M_n$ the family $\{\e^{tA}\mid t\in\mbR\}$ is a one-parameter subgroup of $GL_n,$\footnote{We do not care which field -- $\mbR$ or $\mbC$ -- $GL_n$ is defined over.}  the general linear group (of invertible matrices) of degree $n.$ In other words, the solution of 
\[
\frac{du}{dt} = Au
\]
lives in the Lie group whose Lie algebra is $M_n.$ 
So the local error of the Lie-Trotter and Strang splittings can be written in the terms of Lie brackets (\eqref{eq:2.add.1}, \eqref{eq:2.add.2}) and the splitting is exact if the vector fields of the subsystems commute. 
Note that $\e^{A_1}$ and $\e^{A_2}$ both live in the same Lie group, which is a crucial moment: the flow of a vector field 
$A\in C^\infty$ defines its own group of diffeomorphisms in the general case of non-constant coefficients. Alternatively, we can say $A_1$ and $A_2$ belong to the same Lie algebra.
\end{example}

This algebraic interpretation  goes beyond this finite-dimensional case and is the basis for splitting methods as geometric integrators. 
References are \cite{HaiLu10geometric, BlaCa16concise, McLachQuis02splitting, HundsVer07Numerical, Sanz97geometric, SanzCal94numerical}.

Geometric integration can be defined as developing integrators that preserve some underlying geometric properties and structures such as symmetries, solutions living on a special manifold, first integrals etc. An important example is symplectic methods that preserve, in particular, phase volume\footnote{They also typically preserve important first integrals and  ''almost'' preserve energy in practice~\cite{Sanz97geometric}.} 

We need some basic facts about Hamiltonian systems.

\begin{example}(e.g.~\cite{Hall13quantum, BlaCa16concise, SanzCal94numerical})
Let $H$ be a smooth autonomous Hamiltonian and let $x=(p,q)$ be an element of the corresponding phase space $\mbR^{2m},$ where $p$ are coordinates and $q$ are  momenta. The Hamiltonian equations are
\begin{align*}
\dot{x} &= J \grad H(x), \\
J &=  \begin{pmatrix}
0 & \Id \\ 
-\Id & 0
\end{pmatrix}.
\end{align*}
 The Lie algebra of Hamiltonian functions is the linear space of all smooth functions  endowed with the Poisson bracket as the Lie bracket. 
Any Hamiltonian function defines a Hamiltonian  vector field 
\[
X_f = J \grad f.
\]
Hamiltonian vector fields then form a Lie algebra 
 with the Lie bracket
\[
[X_f, X_g] = X_{\{g,f\}}.
\]
The flow $\Phi_f$ of a Hamiltonian vector field is the solution of  
\[
\dot{\Phi}_f = X_f(\Phi_f).
\]
The solution operator is the corresponding exponential mapping. We always have
\[
(D\Phi_f(x))^T J D\Phi_f(x) = J,
\]
where $D\Phi_f$ is the Jacobian of $\Phi_f.$ 
\end{example}

Transformations having the last property are called symplectic. They form a Lie group\footnote{In  the general case, the Lie group of transformations that preserve a symplectic $2-$form on the symplectic manifold of a physical system.}. In particular, the phase volume is preserved under the composition of Hamiltonian flows and symplectic integrators provided they all share the same phase space, since the composition of symplectic integrators is again symplectic. 

\begin{example}
The one step of the St\"ormet-Verlet/leapfrog splitting is 
\begin{align*}
\wt{p}_{n+1} &= p_n + \frac{h}{2} \frac{\pt H}{\pt q} (\wt{p}_{n+1}, q_n), \\
q_{n+1} &= q_n - \frac{h}{2}\left( \frac{\pt H}{\pt p}(\wt{p}_{n+1}, q_n) + \frac{\pt H}{\pt p}(\wt{p}_{n+1}, q_{n+1}) \right), \\
p_{n+1} &= \wt{p}_n + \frac{h}{2} \frac{\pt H}{\pt q}(\wt{p}_{n+1}, q_{n+1}).
\end{align*}
The scheme is explicit for separable Hamiltonians.
\end{example}

We also need the following  Baker-Campbell-Hausdorff  formula~\cite{BonFul12theorem, BlaCa16concise, Ross02introduction, SanzCal94numerical}\footnote{It is also used in the classical paper~\cite{Ku80representation}.}.  Given elements $a_1, a_2$ of a free Lie algebra the equation in the corresponding Lie group
\begin{equation}
\label{eq:12.1}
\e^{a_1} \e^{a_2} =\e^{b}
\end{equation}
has the solution $b$ that admits a representation as the infinite series
\begin{align*}
b &= a_1 + a_2 + \frac{1}{2}[a_1, a_2] + \frac{1}{12}\left( [a_1, [a_1, a_2]] + [a_2, [a_2, a_1]]\right) + \ldots,
\end{align*}
where all other terms can also be written as combinations of iterated Lie brackets. This expression should be understood as a formal power series in the corresponding algebra of series~\cite{BonFul12theorem, BlaCa16concise}. Such an abstract interpretation is  studied in~\cite{AuHerKochThal17BCH_formula}  in the setting of operator splitting methods.  

Whether this representation makes sense in a particular setting is actually a non-trivial question. In the case of a finite dimensional Lie algebra $a_1$ and $a_2$ should be in a neighborhood of $0,$ and the series can diverge even for matrices~\cite{BlaCa04convergence}. 
Actually, the equation~\eqref{eq:12.1} is only guaranteed to be always solvable in the algebra of formal polynomials  and may not admit the series representation in the terms of commutators of the initial operators etc~\cite{BiaBonMa20Baker, BlaCa16concise, BlaCa04convergence}\footnote{Such pathologies are sometimes ignored in the literature.}. 

Now the idea is clear. Consider a particular case of~\eqref{eq:5.5}
\begin{equation}
\label{eq:12.3}
\sum_{i} a_i \e^{b_{i,1}hA_1} \cdots \e^{b_{i,m}hA_m},
\end{equation}
where $h$ is the size of the time step,  $\sum_{j=\ov{1,m}} A_i = A$ and $A$ is the original operator of a system. Then the BCH formula allows us to find the series $B_i$ such that
\begin{equation*}
\label{eq:12.2}
\e^{b_{i,1}A_1} \cdots \e^{b_{i,m}A_m} = \e^{B_i},
\end{equation*}
where 
\[
B_i= hB_{i,1} + h^2 B_{i,2}+\cdots
\]
 and every $B_{i,k}, k\ge 1,$ is a linear combination of commutators. Expanding every product of exponentials in powers of $h$ in~\eqref{eq:12.3} gives
\begin{align*}
\e^{hA} - \sum_{i} a_i \e^{b_{i,1}hA_1} \cdots \e^{b_{i,m}hA_m} = hA + h \sum_{k} P_{1,k} C_{1,k} + h^2 \sum_{k} P_{2,k} C_{2,k} + \cdots,
\end{align*}
where $P_{i,k}$ are polynomials in $\{a_i\}$ and $\{b_{i,k}\}$ and $\{C_{1,k}\}$ are iterated commutators of $A_1,\ldots, A_m$ ($C_{1,k}$ is simply $A_k$). 
Then the scheme is first order accurate if
\begin{align*}
P_{1,k} &= -1, \quad P_{2,k} = 0, 
\end{align*}
and so on. 

To obtain schemes of the third and higher order some time steps $b_{i,k}$ and coefficients $a_{i}$ should be either negative or complex~\cite{HaiLu10geometric, BlaCa16concise}\footnote{See also for a discussion of stability~\cite{HundsVer07Numerical, Schatz02numerical,CasCharDesVil09splitting}.}. The first condition does not work in the case of diffusion operators (Example~\ref{ex:9.2}), while complex time steps may be actually usable here since the semigroup of a diffusion operator can be analytic (in a cone in $\mbC$). 

\begin{example}(\cite{HundsVer07Numerical})
If $S_h$ denotes one step of the Strang splitting, the scheme whose one step equals
\[
\frac{4}{3} (S_h)^2 -\frac{1}{2}S_h
\]
is (formally) fourth order accurate.
\end{example}

Example~\ref{ex:20} also applies here.

The idea to use the BCH formula and the approach of formal power series in the case of PDEs usually fails~\cite{BlaCa16concise}  as the corresponding smooth algebraic structure is not evident or does not exist,\footnote{This is related to the notion of hypoellipticity and the H\"{o}rmander condition. See also~\cite{BiaBon14convergence, Bon10ode_version}.} and the applicability of the BCH formula is not clear. Numerical approximations are usually discretized, as we have already discussed, and thus do not share geometrical properties of the real solution. The possible presence of BC and incompatible domains which the operators of the system and subsystems are defined on are other problems. Still, one can apply results formally or as the first steps of a rigorous proof (which usually involves bounds for commutators), or simply use the method described above to find new higher order exponential schemes.   

The Lie formalism~\cite{HaiLu10geometric, LanVer99analysis, HundsVer07Numerical, DiFaHaZla01commutativity}\footnote{The Lie formalism can be used to handle the Taylor expansion (e.g.~\cite{DesDuDuLouMa11adaptive}).} is typically used then. That is, given an ODE or a IVP
\[
\frac{du}{dt} = F(u),
\]
the solution admits a formal representation as a Lie-Gr\"{o}bner series
\[
u(t,x) = \e^{tL_F} \Id(x) =\sum_{k\in \mbN} \frac{t^k}{k!} L^k_F \Id(x),
\]
where $L_F$  is the Lie-Gr\"{o}bner operator that acts on a operator $g$ as
\[
L_F g = g^\prime F,
\]
where $g^\prime$ is the Frechet derivative.  $L_F$ is the Lie derivative $\sum_{i} F_i \frac{\pt }{\pt x_i}$ if  $F$ is a vector field. 

Now any exponential splitting schemes can be treated in the terms of such exponentials. In particular, we obtain Examples~\ref{ex:10.1} and \ref{ex:10.2} since the commutator $[L_{f_1}, L_{f_2}]$ is calculated at $\Id,$ so we use $(L_{f_1}L_{f_2} -L_{f_2}L_{f_1}) \Id$.


\begin{example}(\cite{HundsVer07Numerical})
Consider 
\[
\frac{du}{dt} = \tr(A\Hess u) + f(u),
\]
where $f$ is a function. Separating diffusion and reaction via the Lie-Trotter splitting is second order accurate if
\[
\tr(A\Hess f(u)) = f^\prime(u) \tr(A\Hess u).
\] 
This is rarely the case in practice.
\end{example}

A general framework for unbounded generators of $C_0-$semigroups is developed in~\cite{HanOs09exponential}. Other references include \cite{CsoFaHa05weighted, FaGei07iterative, BlaCalCaSanz21symmetrically, Yo90construction, HaSu05finding, Su91general, Su92general,Su95hybrid,Su93improved,Su91general, BlaCaMu08splitting, BlaCaCharMu13optimized, KochLu08variational,CasCharDesVil09splitting}. Sources of the physical nature in Section~8 are also often valid here.

Algebraic calculation around the BCH formula and error expansions of Section~8 can be combined (e.g.~\cite{KochLu08variational}).

However, another approaches to constructing exponential splittings exist. 

\begin{example}
\cite{GyKry05accelerated} develops higher order schemes by combining simpler operator splitting methods on  meshes with different sizes and with variable time steps\footnote{which does not pair well with the BCH formula, for instance} for possibly degenerate second order parabolic equations in Sobolev spaces. 
\end{example}



\begin{remark}
The backward error analysis for exponential splitting methods is given in~\cite{HaiLu10geometric}. This universal approach represents  a numerical scheme as the flow of some equation and uses a truncated version of this flow to study the global error. 
\end{remark}

\begin{remark}
So-called B-series can be used  instead of the BCH formalism~\cite{HaiLu10geometric, BlaCa16concise,AuHerKochThal17BCH_formula}. 
\end{remark}

\section{The It\^o-Taylor expansion, the Kunita expansion and exponential methods for SDEs}

Starting from this point, we only briefly touch the nature of the results cited and mostly concentrate on providing enough references an interested reader could benefit from.

Roughly speaking, classical numerical methods for strong approximations  of SDEs are typically based on iterating the It\^o formula and the corresponding It\^o-Taylor (Wagner-Platen) expansion (alternatively, the Stratonovich-Taylor expansion), and weak approximations  add parabolic equations, the Girsanov theorem and the Malliavin calculus into the picture.  However, the Malliavin calculus and the theory of Lie algebras are extensively involved in the well-established stochastic algebraic theory behind  the Stratonovich-Taylor expansion and its truncated versions~\cite{Ku80representation,  TaKo09operator, Mi01Lie, Ku01approximation, Ku04approximation, Fu06sixth} (also e.g.~\cite{MaNohSil20algebraic, NoVic08weak, Ya18weak, FourLasLeLions99applications, OshiTeichVe12new, Ya17higher, IguYa21operator})\footnote{The number of publications devoted to the algebraic interpretation of the It\^o-Taylor expansion and the corresponding advances in the Malliavin calculus is huge  so all sources listed in this section should merely be seen as starting points or semirandom examples.} in general and in obtaining cubature formulas for weak approximations~\cite{LyonsVic04cubature, CriMaNee13cubature, CriOr13Kusuoka, GyurLyons11efficient, HaTa22MC_construction} in particular, which provides an appropriate counterpart of the theory in the previous section.

We briefly recall the formal starting point. Consider smooth $\mbR^m-$valued $\sigma_0, \sigma_k, k=\ov{1,m},$ with bounded derivatives and the corresponding Lie derivatives
\[
X_k = \sum_{i=\ov{1,m}}\sigma_{k,i}\frac{\pt}{\pt x_i}, \quad i =\ov{0,m}.
\]
Set $w_0(t) =t.$ Then the solution of
\[
d\xi(t) = \sigma_0(\xi(t)) dt + \sum_{j=\ov{1,m}} \sigma_j(\xi(t)) \circ dw_j(t),
 \]
can be represented as
\begin{align*}
\xi(t) &= \e^{\eta(t)}, \\
\eta(t) &= \sum_{i=\ov{0,m}} w_i(t)X_i + \sum_{\alpha \in A} c_\alpha I_\alpha X^\alpha,
\end{align*}
where $A$ is the set of special indexes, $\{c_\alpha\}$ are numbers,  $\{I_\alpha\}$ are iterated Stratonovich integrals w.r.t. $w_0, \ldots, w_m$ and $\{X^\alpha\}$ are iterated commutators of $X_0, \ldots, X_m.$ By a classical result of H.~Kunita, the series is convergent when the Lie algebra $L$ generated by $X_0,\ldots, X_m$ is nilponent. On the other hand, this representation can formally be  understood as the exponential mapping on some non-commutative Lie algebra (say, the algebra of smooth operators). In practice, higher order commutators are discarded and  a truncated version of $\eta$  is considered. 

Another approach uses truncated versions of the Stratonovich-Taylor expansion (roughly speaking, this produces so called Kusuoka approximations, in particular\footnote{Cubature formulas and so called Kusuoka approximations are well known to be closed related internally (cf. \cite{Ku01approximation} and \cite{LyonsVic04cubature}, for instance).}). Since it is hard to calculate such iterated integrals efficiently, cubature formulas replace mathematical expectations of such integrals with  iterated integrals w.r.t a family of functions of bounded variation. 
To find cubature formulas explicitly, Lie algebras are typically used, too. 

So algebraic calculations are one important component needed to construct such weak approximations. For instance, one can still construct exponential splitting schemes. However, the nature of calculations in general can be rather different unless one indeed combines simpler schemes (such as the Ninomiya-Victoir scheme): though exponential mappings that correspond to small time cubature formulae or approximate schemes appear, the series are typically truncated so algebraic calculations are often performed in the terms of  finite sums over some symbolic algebras (cf.~\cite{AuHerKochThal17BCH_formula} where such truncations are also studied).  

Lower order schemes may not even need the Malliavin calculus or symbolic calculations if a test function is smooth and has bounded derivatives~\cite{NoVic08weak}. 

\cite{Mi01Lie, OshiTeichVe12new,TaKo09operator, Fu06sixth, Mi00numerical} directly discuss what can be seen as stochastic composition and exponential splitting methods. \cite{FosReisStrange23highArx} provides a different perspective in the spirit of the rough path theory.

\begin{example}(\cite{TaKo09operator,NoVic08weak})
The well known Ninomiya-Victoir scheme  can actually be seen as an exponential splitting scheme. Recall that for
\[
dx(t) = a(x(t))dt + \sum_{j=\ov{1,m}} \sigma_j(x(t)) \circ dw_j(t)
\]
(a  version of) the  Ninomiya-Victoir scheme is
\[
\ve \e^{ha} \e^{\Delta w_1(h) \sigma_1} \cdots \e^{\Delta w_m(h) \sigma_m} + (1-\ve) \e^{ha} \e^{\Delta w_m(h) \sigma_m} \cdots \e^{\Delta w_1(h) \sigma_1}
\]
where $\ve$ is a Bernoulli random variable and the exponentials $\{\e^{t \sigma_k}\}$ denote flows of ODEs
\[
\frac{du}{dt} = \sigma_k(u), \quad k=\ov{1,m},
\]
 One can read this as expressions as a strong scheme or as a weak scheme (then $\ve=\frac{1}{2}$). 
\end{example}


\section{General results for SPDEs and SDEs}
We start with recalling the setting of the filtration theory, as many advances in operator splitting methods for SPDEs start here and it can be considered as a motivating example. 

\begin{example}
A standard version of the nonlinear filtration problem is as follows. Consider the state $x$ and the noisy measurement $y$ given via
\begin{align}
\label{eq:14.1}
dx(t) &= b(x(t))dt + \sigma(x(t))  dw(t), \notag \\
dy(t) &= h(x(t))dt + dv(t),
\end{align}
where $v$ and $w$ are independent Wiener processes. The unnormalized conditional density $\rho$ of $x$ given $y$ satisfies the Zakai equation
\[
d\rho(t) =L^*(t)\rho(t)dt +  \rho(t) h \cdot dy(t),
\]
where $L^*$ is the adjoint operator of the generator of $x.$ 
\end{example}

We are going to list three fundamental and general results for the Lie-Trotter splitting for SPDEs. 

\begin{example}\label{ex:14.1}(\cite{BenGloRas92Approximations}; see also \cite{BenGloRas90approximation} for a direct proof for the Zakai equation; cf.~\cite{Te68stabilite, EisenHan22variational}) Consider 
\[
dx(t) = A(t,x(t))dt +  B(t, x(t)) dW(t),
\] 
where $W$ is a cylindrical Wiener process and $A,B$ satisfy standard assumptions:  hemicontinuity, monotonicity, coercivity for the deterministic part, Lipschitz continuity for the stochastic part, measurability and controlled growth for either etc. To define a splitting scheme, two equations in some Hilbert space $H$
\begin{align*}
du(t) &= A(t,u(t))dt, \\ 
dy(t) &=  B(t, y(t)) dW(t)
\end{align*}
 are considered in the variational formulation (via the theory of Gelfand triples), and the actions of their propagators are iterated so the deterministic (a PDE) and stochastic (a Markov semigroup) parts of the original equation are combined exactly as in the Lie-Trotter splitting. Define $t_k = \tfrac{k}{n}, k=\ov{1,n},$ for some fixed natural $n$ and let $d_*(t)$ ($d^*(t)$) be the smallest (largest) $t_k$ such that $t_k \le t$ ($t_k>t$). Then consider\footnote{Schemes that swap propagators or combine them on the second stage are also considered.}
\begin{align}
\label{ex:14.5}
u_n(t) = x(0) + \int_0^t A(s, u_n(s)) ds + \int_0^{d_*(t)}\!\! B(s, y_n(s)) dW(s), \notag \\
y_n(t) = x(0) + \int_0^{d^*(t)}\!\! A(s, u_n(s)) ds + \int_0^{t} B(s, y_n(s)) dW(s),
\end{align}
which provides a compressed and convenient representation of the Lie-Trotter splitting on the whole time interval.
Then $u_n, y_n\to x, n\to\infty,$ in $L^2((0;T)\times \Omega; H)$ and $u_n(t), y_n(t)\to x(t), n\to\infty,$ in square mean.  See~Example~\ref{ex:16.1} for ordinary SDEs.  
\end{example}


\begin{example}\label{ex:14.2}
\cite{GonKo98Fractional}  discusses a variety of SPDEs including those with non-classical assumptions and applications of the two step Lie-Trotter splitting to such evolution equations. In contrast to~Example~\ref{ex:14.1}, also schemes for which both steps contain stochastic integrals (w.r.t. to orthogonal noises) are considered. Splitting schemes are used to derive a constructive proof of the existence of the initial SPDE, in particular. Splitting is also used to study invariant sets and to obtain comparison theorems (and then establish conditions for positivity of solutions) (see also \cite{Ko92comparison} for such results as corollaries of the classical formulation of the Trotter-Kato formula). Applications include some types of  McKean-Vlasov equations. 
\end{example}

\begin{remark}
Using splitting to prove existence theorems is a useful trick (see some references in~\cite{GonKo98Fractional}; also~\cite{San13splitting_up, DeuSan14convergence, KoNo16time,NeidSteZa20convergence_2,Na95remarks} at least to some extend). 
\end{remark}

Both previous examples use the semigroup theory and variational methods for PDEs. 

\begin{example}\label{ex:14.3}
A different approach based primarily on probabilistic arguments is developed in~\cite{GyKry03Splitting, GyKry03rate} where general second-order nonlinear SPDEs
\[
du(t) = (A(t)u(t) + f(t))dv(t) + (L(t)u(t) + g(t))dW(t) 
\]
under the assumption of stochastic parabolicity  are considered. Here $W$ is a finite-dimensional Wiener process, $v$ is coordinate-wise increasing continuous adapted process,  $f,g$ are weakly continuous functions with values in some Sobolev space and $A, L$ are second order and first order respectively differential operators with sufficiently smooth coefficients.  With stochastic integrals being properly separated from the rest of the equation, convergence for the multistep Lie-Trotter splitting is established. Results about the rate of the convergence for $\E \sup_{t} \|u(t)-u_n(t)\|_H,$ where $H$ is some Sobolev space, are given, in particular.
\end{example}

All three examples use the notion of the weak  solution.

Examples~\ref{ex:14.1}, \ref{ex:14.2} and \ref{ex:14.3} are very versatile and cover an extremely wide range of applications, but they are not the only works on the topic. Other references that deal with the theory of nonlinear filtration and, in some cases, provide bounds for the rate of convergence include~\cite{FlorLeGland91time_discretization, HopWong86Lie_Trotter, LeGland92splitting_up, CriLobbeSal22application, ZhangZouZhang22splitting_up, Na95remarks,DeuSan14convergence, BagAnLar23energy, LiWangYauZhang23solving}.   In particular, \cite{ItoRo00approximation} studies splitting for the Kushner equation. 

Other results in an abstract setting without direct invocation of the filtration theory include \cite{PadSheng19convergence, CoxNeer10convergence, KoNo16time, BeMa23convergenceArx, EisenStill22randomizedArx, KofLeMeOs21splitting}. \cite{BarRock17splitting} considers the Douglas-Rachford splitting in the variational setting and thus gives a rather rare example of additive splitting in this context.

\begin{example}(\cite{FlorLeGland91time_discretization, LeGland92splitting_up})
Consider the following modification of~\eqref{eq:14.1}:  
\begin{align*}
dx(t) &= b(x(t))dt + \sigma(x(t)) dw(t) + \alpha(x(t))dv(t),\\
dy(t) &= h(x(t))dt + dv(t).
\end{align*}
under the same assumptions as earlier. Set $a =\sigma\sigma^*$ and $d=\alpha\alpha^*$ and define operators
\begin{align*}
L_0 &= \frac{1}{2}\sum_{i,j} a_{ij}\frac{\pt^2}{\pt x_i \pt x_j} + \sum_{i} b_i \frac{\pt}{\pt x_i}, \\
L_1 &= \frac{1}{2}\sum_{i,j} d_{ij}\frac{\pt^2}{\pt x_i \pt x_j}, \\
B &= h + \alpha \cdot \grad.
\end{align*}
Besides a splitting scheme of the initial SPDE defined via the composition of
\begin{align*}
\frac{\pt u}{\pt t} &= L_0^* u, \\
du &= L_1^* udt+ B^* u dy,
\end{align*}
 one can consider splitting in the form of the composition of the associated SDEs as follows. The first SDE is given via the generator $L_0,$ and the second one is related to 
\[
d\xi(t) = \alpha(\xi(t)) (dy(t) -h(\xi(t))dt)
\]
 via the change of measure with a known density.
\end{example}

\begin{remark}
\cite{FlorLeGland91time_discretization} also contains a discussion of the practical availability of such schemes, which is a rather typical additional problem in practice.  
\end{remark}

Mild solutions of SPDEs and the Trotter-Kato formulae are studied in~\cite{Go06Trotter, Go15Trotter, Go20weak, Go21Trotter} (see also references therein).

However, applied problem still provide plenty of examples that are not necessarily covered by the previous universal approaches.

\begin{example}
The stochastic Allen-Cahn equation is 
\[
du = \Big(\frac{1}{2}\Delta u + (u - u^3)\Big) dt + dW,
\]
where $W$ is a cylindrical Wiener process or a $Q-$Wiener process. Here the nonlinear part $u-u^3$ only satisfies the one-sided Lipschitz condition and grows faster than linearly at infinity\footnote{The one-sided Lipschitz condition is an example of a non-Lipshitz assumption for SDEs that usually guarantees stability of numerical schemes and existence of moments of the solution.}. 
\end{example}

\begin{remark}
Another famous example of models that does not satisfy standard assumptions is the FitzHugh-Nagumo model (known in  SDE and SPDE formulations). Moreover, the second-order operator of this model is not parabolic. 
\end{remark}

\begin{example}(\cite{MaStuHig02ergodicity})
Consider 
\[
d\xi(t) = a(\xi(t)) + \sigma(\xi(t)) w(t).
\]
The split-step backward Euler method with step size $\tfrac{1}{n}$ is
\begin{align*}
\wt{\xi}_{k+1} &= \xi_{k+1} + a(\wt{\xi}_{k+1}) \frac{1}{n}, \\
\xi_{k+1} &= \wt{\xi}_{k+1} + \sigma (w\left(\tfrac{k+1}{n}\right) - w\left(\tfrac{k}{n}\right)), \\
\xi_{k+1} &\approx \xi\left(\tfrac{k+1}{n}\right).  
\end{align*}
In particular, this method can be used to preserve ergodic properties of the initial equation even for $a\not\in Lip(\mbR^m)$ (see also \cite{BuckSamTam22Splitting, CheMelTu21diffusion, BreCo22strong, AbVilZy15long, CuiHongSun21structureArx, AbBuck16splitting, MilTre03quasi, AgaMaMe23random} for results about structure-preservation properties of splitting methods in the stochastic setting, including those about invariant measures and symplectic properties and e.g.~\cite{HanKraOs12second} for a deterministic result outside the theory of symplectic integrators). 
\end{example}

\begin{example}\label{ex:14:5}(\cite{BreCuiHong19strong, BreGou19analysis, BreGou20weak})
The following splitting scheme for the stochastic Allen-Cahn equation on the interval $[0;T]$ is proposed:
\begin{align*}
\wt{u}_{n} &= \phi(u_n), \\
u_{n+1} &= S_{T/n} \wt{u}_n + \int_{kT/n}^{(k+1)T/n} S_{(k+1)T/n - s} dW(s),
\end{align*}
where $(S_t)$ is the semigroup of the Laplacian in a proper function space and 
\begin{equation*}
\label{eq:14.2}
\frac{d\phi(t)}{dt} = \phi(t) -\phi(t)^3.
\end{equation*}
Higher rates of convergence require more regularity of the noise.
\end{example}

 The solution considered in the last example and related publications is often a mild solution so the theory of semigroups is heavily involved.

\begin{remark}
The splitting scheme of Example~\ref{ex:14:5} is the Euler-Maruyama scheme of some auxiliary equation. Such an idea appears also in~\cite{ZhaoJu08convergence, BreCoGior24splitting, KoLarLind18discretization}.
\end{remark}


Some other works on the split-step backward Euler method and similar schemes for SPDEs and SDEs are \cite{BuckSamTam22Splitting, ZhaoJu08convergence, BreCoGior24splitting, BreCo22strong, DoerTeich10semigroup, FuKoLarLind18strong, KoLarLind18discretization, CuiHongLiuZhou19strong, CuiHong19strong_weak, WangLuLiu08two_stage, WangLi10split_step, AnCriOtto23uniformArx, GerJourcCle16Ninomiya, IguYa21operator, KuhNeer04Lie,DoerTeich10semigroup, Lejay18Girsanov},\footnote{not necessarily with problematic coefficients and in a wide range of situations} and \cite{FuKoLarLind18strong, CuiHong19strong_weak, WangMeiLiBu19split_step, BreCoUlan23analysisArx} combine splitting   with space discretization. 

Some applications, including those to particle systems, Heston and Langevin models and physics and finance in general, are \cite{CheMelTu21diffusion, ChenReis22flexible, AbVilZy15long, CuiHongSun21structureArx, AbBuck16splitting, Gei17iterative, Bou17CayleyArx, MilTre03quasi, SbaiJour13high, LordWang22convergenceArx, LenkMa15second_order, LenkMa15weak, MaMon22weak, ChenMar21convergenceArx, AdamsDuReis22operator, KeLord23adaptive, BerDobMon23piecewiseArx, DesDuDuLau14analysis, TeLi20strongArx, DoHan14high, NadZa13approximation, BasTah12strong, Shar03splitting, HoutToi16application,San13splitting_up}\footnote{This occasionally includes exponential splitting such as the Ninomiya-Victoir scheme.}.

Operator splitting methods are well known in the rough path theory. We refer to~\cite{HolKarlKenLie10Splitting} for a general survey and literature (also \cite[Chapter 15]{GlowOshYin17Splitting}, \cite{FrizOber11splitting, He23convergence}, both for RDEs and RPDEs). \cite{FosReisStrange23highArx} uses a related approach.

Unfortunately, statistical applications of operator splitting methods 
seem to be mostly undeveloped. See~\cite{PiSamDit22parameterArx, GerJourCle18asymptotics} for some results and references.

\section{A few examples: diffeomorphic flows of SDEs, Brownian web and non-homeomorphic Harris flows}
We end with some rather simple illustrations in the terms of flows of ordinary SDEs, that is, semigroups $(\phi_{s,t})_{s\le t}$ of random transformations of $\mbR^m$ generated by
\[
d\phi_{s,t}(x) = a(\phi_{s,t}(x))dt + \sum_{j=\ov{1,d}} \sigma_j(\phi_{s,t}(x)) dw_j(t),
\]
where $w_j, j=\ov{1,d}$ are $\mbR^m-$valued Wiener processes.  

The Lie-Trotter splitting that separates drift and diffusion is given in~\eqref{ex:14.5}.
\begin{remark}
The It\^o formula cannot be applied to the processes constructed via this scheme, in contrast to, say, the Euler-Maruyama scheme. 
A method to overcome this difficulty is proposed in~\cite{GyKry03Splitting} where $y_n$ and $u_n$ are combined to introduce an approximate process that satisfies some SPDE. For that, given $u_n, y_n$ on the uniform partition of $[0;T]$ define a process $z_n$ on $[0;2T]$ via
\begin{align*}
z_n(t) &= u_n(\tau_n(t)) + y_n(\kappa_n(t)), \\ 
\tau_n(t) &= \Big(\frac{t}{2}-\frac{k}{n}\Big) \1\Big[ t\in \Big[\frac{2k}{n}; \frac{2k+1}{n}\Big)\Big] + \frac{k+1}{n} \1\Big[ t\in \Big[\frac{2k+1}{n}; \frac{2k+2}{n}\Big)\Big], \\
\kappa_n(t) &= \tau_n\Big(t-\frac{1}{n}\Big).
\end{align*} 
Then
\[
dz_n(t) = a(z_n(t)) d\tau_n(t) + \sum_{j=\ov{1,d}} \sigma_j(z_n(t)) \cdot dw_j(\tau_n(t)).
\]
\end{remark}

\begin{example}\label{ex:16.1}
If $a, \sigma_j, j=\ov{1,m},$ are Lipschitz continuous, the SPDE setting in~\cite{BenGloRas92Approximations} applies, and the rate of convergence is known:
\begin{align*}
\sup_{x\in [-M,M]} &\E \sup_{t\in [0,T]} \left( y^n_t(x) -\phi_{0,t}(x)\right)^2 \le C \delta_n, \\
\sup_{x\in [-M,M]} &\sup_{t\in [0,T]} \E\left( u^n_t(x) - \phi_{0,t}(x)\right)^2 \le C \delta_n,
\end{align*}
\end{example}

Actually, similarities between the Euler-Maruyama scheme and the Lie-Trotter splitting allow one to transfer some proofs for the former scheme to the latter with minor corrections.

\begin{example}
If we consider the one-dimensional case and  Holder continuous coefficients, the version of the Yamada-Watanabe method  developed in~\cite{GyRa11Note} can be used to estimate the rate of the convergence. 
\end{example}

\begin{example}\footnote{This observation seems to be new.}
We can use the classical proof of the Talay-Tubaro expansion~\cite{GraTa13stochastic}, the only difference being that we need to use pathwise Taylor expansions instead of the It\^o formula. E.g. in the one dimensional case, this yields, for smooth $f$ with bounded derivatives,
\begin{align*}
\E f(\xi(t)) - \E f(y_n(t)) = \frac{t}{n} \E \int_0^t \psi(s, \xi(s)) ds + o\Big(\frac{1}{n}\Big)
\end{align*}
where
\[
\psi = \frac{1}{2} \left[ a^2 \pt^2_{xx} +  a a^\prime \pt_x + a\pt_t \pt_x\right] u + \frac{1}{4} \sigma^2 \pt^2_{xx} \left( a\pt_x u \right).
\]
\end{example}

The Lie-Trotter splitting can also applied to non-homeomorphic flows, in which case the corresponding random transformations can be even piecewise constant functions. 

\begin{example}(\cite{DoVov18Arratia, DorVov20ApproximationsEng})
An extreme example is the Brownian web $\{B_{s,t}\},$ a system of Wiener processes in $\mbR$ that start from all points of $\mbR$ at all times, move independently until they meet and merge afterward. Here $B_{s,\cdot}(x)$ is a Wiener process started at $x$ at time~$s.$ Though the Brownian web does not admit a representation as a SDE, the Lie-Trotter splitting can still be defined on the interval  $[\frac{k}{n}; \frac{k+1}{n})$ as
\begin{align}
\label{eq:16.1}
u_n(t) &= S_{t-k/n} \, y_{n-1}\left(\tfrac{k}{n}\right), \notag \\
y_n(t) &=  \left( B_{t, k/n} \circ S_{1/n} \right) u_n\left(\tfrac{k+1}{n}\right),
\end{align} 
where $\frac{dS_t}{dt}=a(S_t).$ Then the flows $\{y_n\}$ converge in distribution to the Brownian web with drift $a$ and so do the associated pushforward measures. 
\end{example}

\begin{example}(\cite{Vo22SplittingManu})
A Harris flow $X_{s,t}$ is a flow of one-dimensional correlated Wiener processes  in which 
\[
\frac{d}{dt}\langle X_{s,\cdot}(x_1), X_{s,\cdot}(x_2)\rangle(t) = \phi\left(X_{s,t}(x_1)-X_{s,t}(x_2)\right),
\]
for some symmetric strictly positive definite $\phi\in C(\mbR).$  Under proper assumptions $X_{s,t}, t > s,$ is a piecewise constant function. Moreover, though $X_{s,\cdot}(x)$  admits a representation as a solution to some SDE, the process $X_{s,\cdot}(x)$ may not be the unique solution~\cite{WarWa04Spectra}. Still, one defines the Lie-Trotter splitting as in~\eqref{eq:16.1} and the resulting flows converge in distribution to the flow with the same $\phi$ and the corresponding drift\footnote{Proper numerical schemes for such non-homeomorphic flows are studied in~\cite{Gli12Disordering, Nish11Discrete}.}. Pushforward measures and dual flows in the reversed time are also convergent.  
\end{example}


\section*{Acknowledgments} 
The author expresses his deepest gratitude to the Armed Forces of Ukraine for keeping him safe and thus making this work possible.   

This work was supported by a grant from the Simons Foundation (1030291, M.B.V.).
\printbibliography

\end{document}